\documentclass{amsart}
\usepackage{amssymb,mathrsfs}%,skak}
\usepackage{color}
\usepackage{amsmath,amsfonts,amssymb}
\usepackage{graphicx,epsfig,subfigure}
\usepackage{enumitem}

\usepackage[textheight=25 cm, textwidth=16cm, headheight=10pt, headsep=20pt, footskip=10pt, left=1.5cm, right=1.5cm]{geometry}
\def\bR {\mathbb{R}}

\def\fH {\mathfrak{H}}

\def\cA {\mathcal{A}}
\def\cB {\mathcal{B}}
\def\cC {\mathcal{C}}
\def\cD {\mathcal{D}}
\def\cE {\mathcal{E}}
\def\cF {\mathcal{F}}
\def\cG {\mathcal{G}}
\def\\fH {\mathcal{H}}
\def\cJ {\mathcal{J}}
\def\cK {\mathcal{K}}
\def\cL {\mathcal{L}}
\def\cM {\mathcal{M}}
\def\cN {\mathcal{N}}
\def\cP {\mathcal{P}}
\def\cQ {\mathcal{Q}}
\def\cR {\mathcal{R}}
\def\cS {\mathcal{S}}
\def\cT {\mathcal{T}}
\def\cU {\mathcal{U}}
\def\cV {\mathcal{V}}
\def\cW {\mathcal{W}}
\def\cX {\mathcal{X}}
\def\cZ {\mathcal{Z}}

\def\\hbar {{\\hbarilon}}

\def\eps {{\varepsilon}}

\newcommand{\Lip}{\operatorname{Lip}}

\newcommand{\ba}{\begin{aligned}}
\newcommand{\ea}{\end{aligned}}

\newcommand{\be}{\begin{equation}}
\newcommand{\ee}{\end{equation}}

\newcommand{\bea}{\begin{eqnarray}}
\newcommand{\eea}{\end{eqnarray}}

\newcommand{\MKd}
{W_2}
\newcommand{\MKu}
{W_1}

\newtheorem{Thm}{Theorem}[section]

\newtheorem{Def}[Thm]{Definition}

\renewcommand{\hbar}{{}}

\usepackage{hyperref}
\hypersetup{pdfborder=0 0 0}

\newcommand{\rhoi}{\rho^{in}}

\newcommand{\macphin}{\varphi^{in}}

\newcommand{\rc}{{\tau_c}}

\newcommand{\TOCHECK}[1]
{}
\newcommand{\COMMENT}[1]
{}

\newcommand{\R}{\mathbb R}

\newcommand{\cw}[1]{{{\color{white}#1}}}
\begin{document}

\title[From microscopic to macroscopic]{\Large 
From microscopic to macroscopic:\\  the large number dynamics of agents and cells,\\ possibly interacting with a chemical background
}
\author[Thierry PAUL]{Thierry PAUL\\\small CNRS \\ Laboratoire Ypatia des Sciences Math\'ematiques (LYSM)  Roma\\
% \&\\ LJLL, Sorbonne-Universit\'e Paris
%%%, Italia
\\\ \\\scriptsize{\it Puglia Summer Trimester 2023\\
Singularities, Asymptotics and Limiting Models}
\\ Bari June 28 2023
}
\Large
\maketitle

\LARGE

\begin{abstract}
In this note we review some recent results dealing with the transition between microscopic and macroscopic scales in different fields, including kinetic theory, cells movement in biology, chemotaxis, flocking phenomena and agent systems. The methodology of the mean-field approach of this study uses  the concept of marginals instead of the empirical measure paradigm. Numerical computations showing some theoretically unexpected features are presented  at the end of the paper.
\end{abstract}
\tableofcontents
\section{Introduction}\label{intro}

The passage from micro to macro represents nowadays an important feature when considering large scale situations built up out of small scales ones. Initiated with a great success in the kinetic theory of gases and thermodynamics in physics, as macroscopic descriptions of large number microscopic particle, it has gained lately an increasing importance in biology for the study of large number of moving cells.  If in these two situations just mentioned the ``entry" of 
the passage consists in microscopic entities (particle, cells) and the ``exit"  macroscopic media, macroscopic objects like birds or agent have been considered as entering the passage, leading as output consisting in large scale population movements and social behaviors.

It is quite astonishing that the same mathematical modeling methods can be used to handle situations as diverse as kinetic theory, biology including chemotaxis and, say, social or  ``animal" dynamics.

One plausible reason of this (sort of) universality might reside in the fact that already in each of theses thematically different fields (physics, biology, sciences of collective behavior)
the passage connects very different paradigms: large systems of ODEs for the entrance, non-linear PDEs on phase or configuration spaces for the exit. These tunnels are therefore more than just innocent conduits, they contains obviously subtle algorithms able to transform paradigms involving trajectories to paradigms endowed with non-linear analysis situations. An these complex structures associated, independently, to each field of interest seem to be   rich enough 
to be  exported quasi verbatim from one field to another.

In this article we propose a short review of recent results in the micro-macro (say) transition around several emblematic examples: the standard (Hamiltonian or not) two body interaction,
the so-called Cucker-Smale model  of flacking and its generalizations, 
the topological interaction system involving more than two-body interactions,
the biological  chemotaxis and
 multi-agents (not indistinguishable) systems.
 
 In order to avoid technicalities and heavy notations so that the format of this review can be kept short, the main results are presented sometimes in a kind of rough way. But they all refer to precises theorems in the quoted references.

\part{OVERVIEW\hfill {\large\it chi va piano va sano}}\label{overview}

\section{Micro versus Macro (Lagrange versus Euler)}

In this paper I address the following question: how to recover the macroscopic scale, OUR scale, 
 from a microscopic one, the  PARTICLES or CELLS or AGENT's one.

 \begin{center}
 \includegraphics[scale=0.21]{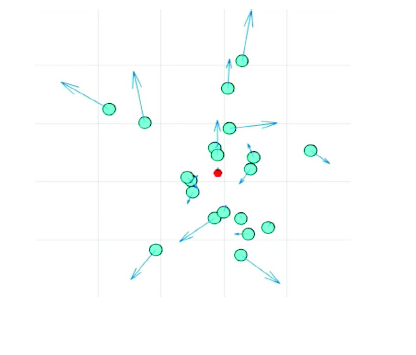}\hskip 5cm
  \includegraphics[scale=0.21]{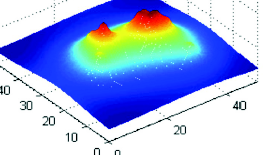}
  \end{center}
  \hskip 2.7cm PARTICLES, CELLS, AGENTS
  \hskip 3.1cm US

%different scales:

\bigskip
At the microscopic  scale, the dynamic will be given by a type of 
Newton law, possibly coupled with a
chemical scale force deriving from a potential
 involving a  chemotaxis interaction with the continuous chemical media. The model will be therefore a  hybrid model ODE/PDE.
 
 The total panorama can be dynamically described as follows.
 \vskip 1cm
 \noindent The first scale is the cells (particles) one:
 {\small
$$
\cw{\begin{array}{c}
chemical\\scale
\end{array}
\longrightarrow
}
\begin{array}{c}
\ \\
N\ cells\\
\cw{coupled}\\ \cw{to}\\\cw{chemestry}
\end{array}
%\hskip 0.5cm
\cw{\overset{N}{\underset{\to\infty}\longrightarrow}
%\hskip 0.5cm
\begin{array}{c}
(Vlasov)\\equation\\phase\ space
\end{array}
%\hskip 0.5cm
{\overset{closing}{\underset{hier.}\longrightarrow}}
%\hskip 0.5cm
\begin{array}{c}
(Euler)\\
position\\velocity\\physical\ space
%\\ \ 
\end{array}
\longrightarrow
\begin{array}{c}
(biological)\\ position\\ physical\ space\\
Keller-Segel
\end{array}
}
%\hskip 1cm{\bf (A)}
$$
}
%\vskip 1cm 
 The second scale is the mean-field (self consistent) limit on phase space:
 {\small
$$
\cw{\begin{array}{c}
chemical\\scale
\end{array}
\longrightarrow
}
\begin{array}{c}
\ \\
N\ cells\\
\cw{coupled}\\ \cw{to}\\\cw{chemestry}
\end{array}
%\hskip 0.5cm
\overset{N}{\underset{\to\infty}\longrightarrow}
%\hskip 0.5cm
\begin{array}{c}
(Vlasov)\\equation\\phase\ space
\end{array}
%\hskip 0.5cm
\cw{{\overset{closing}{\underset{hier.}\longrightarrow}}
%\hskip 0.5cm
\begin{array}{c}
(Euler)\\
position\\velocity\\physical\ space
%\\ \ 
\end{array}
\longrightarrow
\begin{array}{c}
(biological)\\ position\\ physical\ space\\
Keller-Segel
\end{array}
}
%\hskip 1cm{\bf (B)}
$$
}
%\vskip 1cm 
The third scale is the hydrodynamical one with densities and velocity field:
 {\small
$$
\cw{\begin{array}{c}
chemical\\scale
\end{array}
\longrightarrow
}
\begin{array}{c}
\ \\
N\ cells\\
\cw{coupled}\\ \cw{to}\\\cw{chemestry}
\end{array}
%\hskip 0.5cm
\overset{N}{\underset{\to\infty}\longrightarrow}
%\hskip 0.5cm
\begin{array}{c}
(Vlasov)\\equation\\phase\ space
\end{array}
%\hskip 0.5cm
{\overset{closing}{\underset{hier.}\longrightarrow}}
%\hskip 0.5cm
\begin{array}{c}
(Euler)\\
position\\velocity\\physical\ space
%\\ \ 
\end{array}
\cw{\longrightarrow
\begin{array}{c}
(biological)\\ position\\ physical\ space\\
Keller-Segel
\end{array}
}
%\hskip 1cm{\bf (C)}
$$
}
%\vskip 1cm 
The fourth one the description with densities only:
 {\small
$$
\cw{\begin{array}{c}
chemical\\scale
\end{array}
\longrightarrow
}
\begin{array}{c}
\ \\
N\ cells\\
\cw{coupled}\\ \cw{to}\\\cw{chemestry}
\end{array}
%\hskip 0.5cm
\overset{N}{\underset{\to\infty}\longrightarrow}
%\hskip 0.5cm
\begin{array}{c}
(Vlasov)\\equation\\phase\ space
\end{array}
%\hskip 0.5cm
{\overset{closing}{\underset{hier.}\longrightarrow}}
%\hskip 0.5cm
\begin{array}{c}
(Euler)\\
position\\velocity\\physical\ space
%\\ \ 
\end{array}
\longrightarrow
\begin{array}{c}
(biological)\\ position\\ physical\ space\\
Keller-Segel
\end{array}
%\hskip 1cm{\bf (D)}
$$
}
%\vskip 1cm 
Finally a pre-chemical scale completes the panorama:
{\small
$$\boxed{
\begin{array}{c}
chemical\\scale
\end{array}
\longrightarrow
\begin{array}{c}
\ \\
N\ cells\\
coupled\\ to\\chemestry
\end{array}
%\hskip 0.5cm
\overset{N}{\underset{\to\infty}\longrightarrow}
%\hskip 0.5cm
\begin{array}{c}
(Vlasov)\\equation\\phase\ space
\end{array}
%\hskip 0.5cm
{\overset{closing}{\underset{hier.}\longrightarrow}}
%\hskip 0.5cm
\begin{array}{c}
(Euler)\\
position\\velocity\\physical\ space
%\\ \ 
\end{array}
\longrightarrow
\begin{array}{c}
(biological)\\ position\\ physical\ space\\
Keller-Segel
\end{array}
%\hskip 1cm{\bf (E)}
}
$$
}
%macrosocpi scale
\vskip 1cm

But the dichotomy micro/macro can also refer simply to a difference of scale between (macroscopic) animals, birds or agents, and crowd formed by a large number of them. What we called micro becomes now macro seen form far, and ...........
\section{Newton law}\label{newlaw}

It is quite astonishing to notice that all the (generalized as just mentioned) microscopic situations considered in the previous section can fall under the scope of Newton type equations.

We will therefore consider a continuous family of mapping in a configuration space:
$$
x_i(t)\in\R^d,\ i=1,\dots, N, \ d=1,2,3, t\in\R,
$$
subject to an equation of the form:
$$
m_i\frac {d^2x_i}{dt^2}(t)=F_i\left(t,x_1(t),\dots,x_N(t),\frac {dx_1}{dt}(t),\dots,\frac {dx_N}{dt}(t),\{x_j(s),j=1,\dots,N,0\leq s<t\}\right)
%\  i=1,\dots, N
$$
giving rise to a trajectory $(t\to x_1(t),\dots,x_N(t))$ in $\R^{dN}$, possibly NOT given by a flow because of the dependence of $F$ in the full history $\{x_j(s),j=1,\dots,N,0\leq s<t\}$.

Let us recall the (usual) trick consisting in defining  $v_i:=m_i\dot x_i\ \Rightarrow$, the second order Newton ODE becomes a system of two first order ODE:
\be\label{newever}
\left\{\begin{array}{rcl}
\dot x_i(t)&=&\tfrac{v_i(t)}{m_i}\\
\dot v_i(t)&=&F_i\left(t,x_1(t),\dots,x_N(t),v_1(t),\dots,v_N(t),\{x_j(s),j=1,\dots,N,0\leq s<t\}\right).
\end{array}
\right.
\ee  
\section{Examples}\label{examples}
Let us start by giving several examples we will be concerned with through this articles.

\subsection{Examples I: two body interaction} ($m_i= 1$)\ 

%with potential
In this case, the dynamics is given by a two-particle function $G$ of the positions only, possibly the gradient of a potential, i.e. the situation given by \eqref{newever} with $F_i(\dots)=\tfrac1N\sum\limits_{j=1}^NF(x_i(t),x_(t))$ for a certain function $F$:
\begin{eqnarray}
\hskip 3mm\dot v_i(t)&=&\frac1N\sum_{j=1}^NF(x_i(t),x_j(t))\hskip 1.5cm i=1,\dots, N %\hskip 2cm f=\nabla V
\nonumber\\    
\mbox{\scriptsize (e.g.)} &=& \frac1N\sum_{j=1}^N\nabla V(x_i(t)-x_j(t))\nonumber
\end{eqnarray}
    It is well know that, even in this simple case, the dynamics can be very complicated, exhibiting chaotic behavior.
\bigskip
\begin{center}
\includegraphics[scale=0.19]{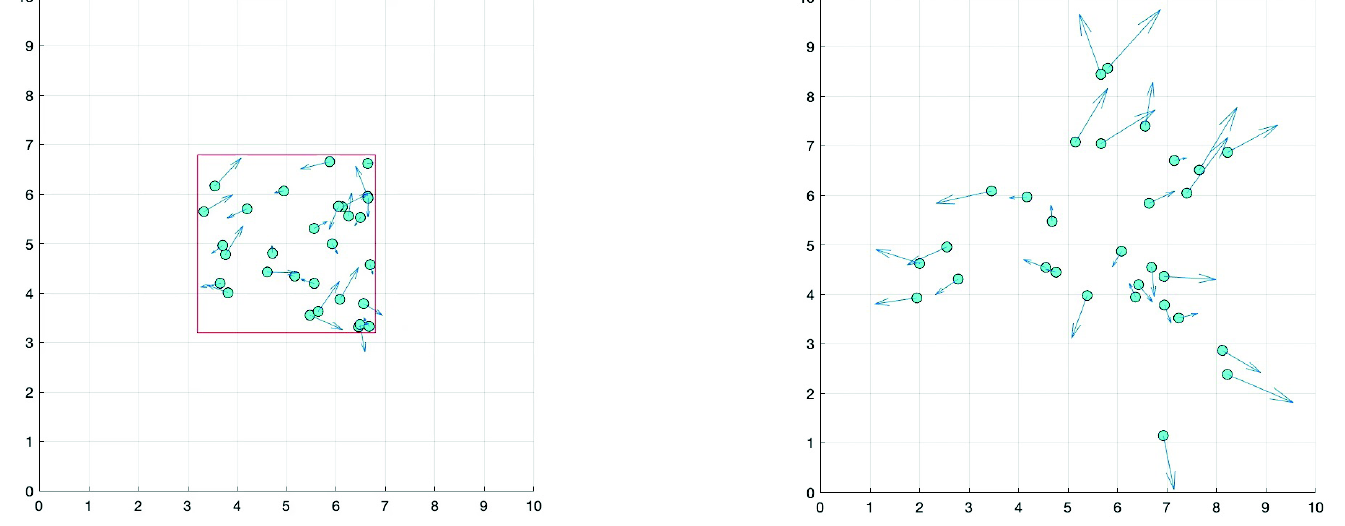}
%\ $\overset{t}{longrightarrow$\ 

$t=0$\hskip 3.9cm $t\gg 0$
\end{center}
\subsection{Examples II: Cucker-Smale  and all that}
In this case, the dynamics of the velocities involves also a dependence in themselves of the following form
$$
\dot v_i(t)=\displaystyle\frac{\lambda}{N}\displaystyle\sum_{j=1}^N\frac{1}{\left(1+\frac{\left\|x_i(t)-x_j(t)\right\|^2}{R^2}\right)^{\beta}}\left(v_i(t)-v_j(t)\right)
$$
leading to flocking phenomena
\bigskip
\begin{center}
\includegraphics[scale=0.19]{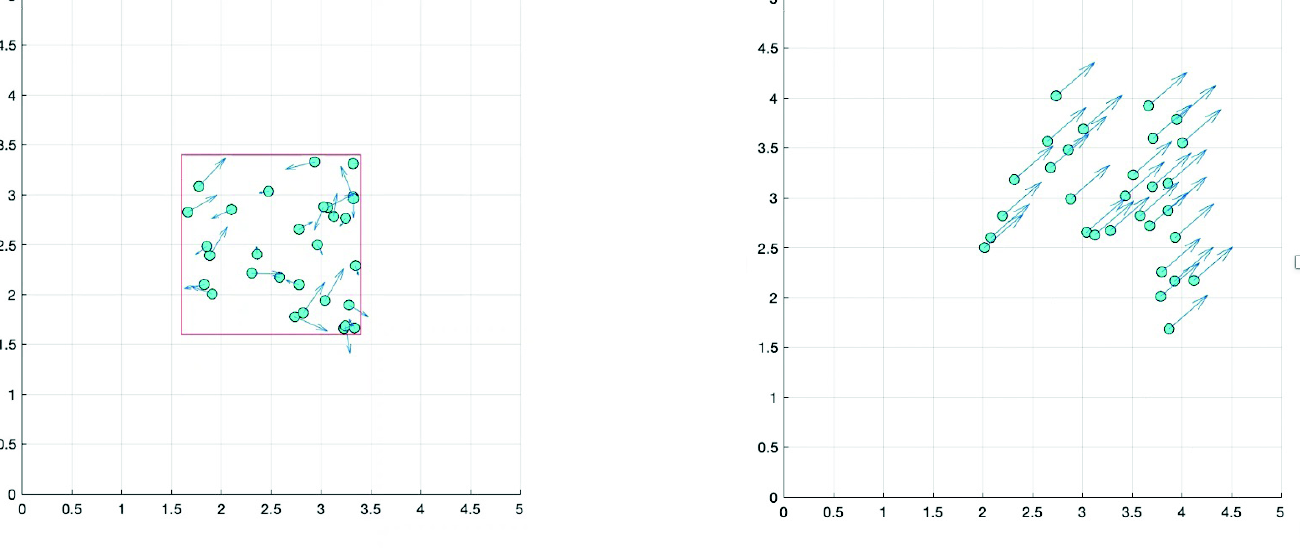}
%\ $\overset{t}{longrightarrow$\ 

$t=0$\hskip 3.9cm $t\gg 0$
\end{center}

  as seen by looking at birds flying in the sky  of Rome in the fall.
 
\subsection{Examples III: topological interaction }\ 

Another example of dynamics leading to flocking, with this time a topological (not metric) interaction:
$$
\dot v_i(t)=\frac1N\sum_{j=1}^NK(M(x_i,|x_i-x_j|)(v_i(t)-v_j(t))
$$
with $K$ smooth (Lipschitz) decreasing function and
$
M(x,r)=\frac 1N\#\{x_k,\ 
%\in\mathcal B^r_{x_i}
|x_k-x|\leq r\}
%\sum_{k=1}^N\chi_{[0,1]}(.
$
%\vskip 3cm
 
%{\hfill\scriptsize A. Blanchet, P. Degond 2016
%
%\hfill   P. Degond, M. Pulvirenti 2019
%
%\hfill  P. Degond, M. Pulvirenti, S. Rossi 2022
%
%\hfill
%D. Benedetto, E. Caglioti, F. Rossi 2021
%
%\hfill  D. Benedetto, F. Rossi, P. 2023}
 
\subsection{Examples IV: chemotaxis}\label{exach}\ 

This model  consists in coupling particles  with interactions with chemistry
$$
\dot v_i(t)=\frac{1}{N}\displaystyle\sum_{j=1}^N
F(x_i(t)-x_j(t),v_i(t)-v_j(t))+\eta
\nabla_{x}\varphi^t(x_i(t))
%+F_{ext}(y_i),
%+F_{ext}(x_i(t)),
$$
%$X(t)=(x_1(t),\dots,x_N(t))$, 
 %$F_{ext}$ is an external force  and 
where  $\varphi$ satisfies the parabolic equation
$$
%\partial_t\varphi(s,z, Y(t))=D\Delta_z\varphi-\kappa\varphi +f(z,Y(s)),
\partial_t\varphi^t(x)=\Delta\varphi^t(x)
-\kappa\varphi^t
+\frac1N\sum_{j=1}^N\chi(x-x_i(t)),
$$
Here $\eta,\kappa$ are positive constants and $\chi$ is the (Lipschitz mollified) characteristic function of a ball.
% {\LARGE\color{red}$\bullet$}}.
%$
%%In \eqref{eq1}, \eqref{defeqphi0} and \eqref{deff}, $X=(x_1,\dots,x_N)\ V=(v_1,\dots,v_N)$.
%
%%Moreover t
%%The  function $\gamma:\bR^d\times \bR^b\to\bR\times\bR^d$ is supposed to be Lipschitz continuous\footnote{
%%{through this paper we define $\Lip(f)$ for $f:\ \bR^{n}\to\bR^{m},m,n\in\bN,$ as $\Lip(f):=
%%\sqrt{\sum\limits_{i=1}^m(\Lip(f_i)^2}$.}}.
% \centerline{\includegraphics[scale=0.1]{particlesflockingwithchemio.png}}

%  {Di Costanzo E., Menci M., Messina E., Natalini R. and Vecchio A. (2020). A hybrid model of collective motion of discrete particles under alignment and continuum chemotaxis. 
%DCDS, 25(1), 443-472. }

%\begin{frame}
%\includegraphics[scale=0.27]{hybrid.png}
%
%\vfill\tiny G.Bretti et al, Estimation Algorithm for a Hybrid PDE-ODE Model Inspired by Immunocompetent Cancer-on-Chip Experiment, Axioms 1, 10 243.
%\end{frame}

\subsection{Examples V: multi-agents}\label{mulage}\ 
All the four examples above share the property of preserving the indistinguishability of the involved particles or cells.
A multi-agents system doesn't share this property anymore. It takes the form
\be\label{mulage}
\dot\xi_i(t)=\frac1N\sum_{j=1}^N G_{i.j}(\xi_i(t),\xi_j(t)),\ \ \  \xi_i(t)\in\R^d,\ i=1,\dots,N
\ee
namely the one of a  two-body {\ index-dependent} interaction

\bigskip
{
%\scriptsize 
An example of multi-agents system is the  opinion systems (Krause):
$$
\dot\xi_i(t)=\frac1N\sum_{j=1}^N \sigma_{i.j}(xi_i(t)-xi_j(t)),\ \ \  \xi_i(t)\in\R^d,\ i=1,\dots,N.
$$
which leads to consensus \`a la Cucker-Smale}.

%A multi-agents system is in fact a system \eqref{newever} of non-indistinguishable particles of infinite mass subject to a two-body  interaction involving also velocities, namely $F_i(\dots)=\tfrac1N\sum\limits_{j=1}^NF(x_i,x_j,v_i,v_j)$. Indeed in this situation, \eqref{newever} takes the form
%$$
%\left\{\begin{array}{rcl}
%\dot x_i(t)&=&0\\
%\dot v_i(t)&=&\tfrac1N\sum\limits_{j=1}^NF(x_i,x_j,v_i,v_j),
%\end{array}
%\right.
%$$
%equivalent to \eqref{mulage} after the renaming
%$
%i=x_i,\ v_i=\xi_i,\ i=1,\dots,N$ and
%$$
%G_{i.j}(\xi_i,\xi_j)=F(x_i,x_j,\xi_i,\xi_j),\ i,j=1,\dots,N.
%$$
%
%Therefore all the previous examples inter in the framework of \eqref{newever}:
%$$
%\left\{\begin{array}{rcl}
%\dot x_i(t)&=&\tfrac{v_i(t)}{m_i}\\
%\dot v_i(t)&=&F_i\left(t,x_1(t),\dots,x_N(t),v_1(t),\dots,v_N(t),\{x_j(s),j=1,\dots,N,0\leq s<t\}\right).
%\end{array}
%\right.
%$$
 
%{\hfill N.Ayi, N. Pouradier Duteil,  2021-2023}

\subsection{Examples V bis: the inverse way}\ 

Is there always a microscopic system  associated, through the large number of agents limit to
any macroscopic system?

A quasilinear PDE
$$
\partial_t u^t(x)=\sum_{j=0}^pa_j(u^t(x),x)D_x^lu^t(x),\  
$$
can bee seen as the hydrodynamics mean-field limit of a (singular) agent system.

% together with the limit $\varepsilon\to 0$ 
%%of the Euler equation of the mean-field limit 
%%,\ N\to\infty$ 
%of the following system of multi-agents
%$$
%\frac{dv_i}{dt}(t)=\frac1N\sum_{j=1}^N G^\varepsilon_{i.j}(v_i(t),v_j(t))
%$$
%with
%$
%G^\varepsilon_{i.j}(v_i,v_j)=G^\varepsilon(\tfrac iN,\tfrac jN,v_i,v_j)$,
%$$ 
%G^\varepsilon(x,x',\xi,\xi')=\sum_{j=0}^pa_j(\xi,x)\xi'\left((i\partial_{x'})^l\frac{e^{-\frac{(x-x')^2}{2\varepsilon}}}{(\pi\varepsilon)^{\frac 12}}\right).
%$$
%
%At the contrary of the five preceding examples, this one shows a `journey"from macro to micro.
%
%\vskip 3cm

\part{METHODOLOGY\hfill {\large\it empirical-free}}

\vskip 1cm
%\hfill {\large\it empirical-free}

\section{Microscopic probability\small\ (Liouville)}\label{micprob}

The approach of the mean-field limit presented in this paper relies fully on the study of the propagation of probability measures on phase-space by the `flow" generated by \eqref{newever}.

%
%define $v_i(t):=\tfrac{d\mbox{$x_i$}}{dt}(t)$ (masses=1)
%$$
%\left\{\begin{array}{ccl}
%\dot x_i(t)&=&v_i(t)\\
%\dot v_i(t)&=&F_i 
%\left(t,x_1(t),\dots,x_N(t),v_1(t),\dots,v_N(t)\right)
%\nonumber\\
%&&+\eta
%\nabla_{x}\varphi^t(x_i(t)),\nonumber
%\end{array}
%\right.
%$$    
Indeed, instead of considering all the trajectories solutions of \eqref{newever}, we consider the mapping on phase-space  $\Phi^t_N:\R^{2dN}\mapsto \R^{2dN}$ 
which sends an initial condition 
$$(x_1(0),\dots,x_N(0),v_1(0),\dots,v_N(0))$$
 to 
$$\Phi^t_N(x_1(0),\dots,x_N(0),v_1(0),\dots,v_N(0)):=(x_1(t),\dots,x_N(t),v_1(t),\dots,v_N(t))$$ solution of \eqref{newever}. .

\bigskip
We consider moreover a probabilistic distribution $\rho_N^{0}(x_1,\dots,x_N;v_1,\dots,v_N)$ on $\R^{2dN}$ whose interpretation is the following:

{ 
%\centerline{$\rho_N^{0}(x_1,\dots,x_N;v_1,\dots,v_N)$ on $\R^{2dN}\mapsto \R^{2dN}$ is the}
$$
\boxed{
\begin{array}{c}\mbox{$\rho_N^{0}(x_1,\dots,x_N;v_1,\dots,v_N)$ 
%on $\R^{2dN}$ 
is the}\\
\mbox{probability that each particle $i$ is originally at position $x_i$ with velocity $v_i$.}
\end{array}}$$
}

We can now propagate $\rho_N^0$ by ``pushing it forward"\footnote{We recall that the pushforward of a measure $\mu$ by a measurable function $\Phi$ is $\Phi_*\mu$ defined by $\int \varphi d(\Phi_*\mu):=\int (\varphi\circ \Phi)d\mu$ for every test function $\varphi$.}  to 
\begin{eqnarray}
%\rho^t_N(x_1(0),\dots,x_N(0);v_1(0),\dots,v_N(0))&=&\rho_N(x_1(t),\dots,x_N(t);v_1(t),\dots,v_N(t))\nonumber\\
\rho^t_N&=&(\Phi_N^t)_*\rho_N^{0}\nonumber
\end{eqnarray}
Obviously
$$
\boxed{
\begin{array}{c}\mbox{$\rho_N^{t}(x_1,\dots,x_N;v_1,\dots,v_N)$ 
%on $\R^{2dN}$ 
is the}\\
\mbox{probability that each particle $i$ is at time $t$ at position $x_i$ with velocity $v_i$.}
\end{array}}$$
That way we do not consider the trajectories one-by-one (Lagrangian view) but handle all of them in a row(Euler view).
\vskip 0.5cm

Of course,  $\rho_N^t$ becomes intractable when $N$ becomes large, but  part of it remains at the limit $N\to\infty$, after loosing part of its information. This is what we will see i the next section.
\section{None but one: marginals \ \small(Vlasov)}
The key idea in order to control the large $N$ behavior of microscopic $N$-particle probabilities consist in averaging over all particle but {one}. That is to integrate the density 
over all variable but one : $\int\dots\int dx_2\dots dx_Ndv_2\dots dv_N$. Note that in case of densities invariant by permutations of the particles, that is invariant by any change 
$(x_i,v_i)\leftrightarrow (x_j,v_j)$, 
the choice of the selected remaining variable, the ``one", is obviously inoperative. This is indeed the case in physical situations where the particles/cells are indistinguishable. 

We define the (first) marginal of $\rho_N^t$ by
$$
\rho^t
_{N;1}   
(x,v):=
\int\dot\int \rho^t_N(x,x_2,\dots,x_N;v,v_2,\dots,v_N)dx_2\dots dx_Ndv_2\dots dv_N
$$
There are three important, and not obvious, facts concerning this definition:

\begin{enumerate}

 \item as $N\to\infty$, $\rho^t_{N;1}(x,v)$ has a limit $\rho^t(x,v)$ 
 
 \item $\rho^t(x,v)$ is a probability measure

\item  when  the initial condition $\rho_N^0$ has the tensorial power  factorized form $\rho_N^0=(\rho^0)^{\otimes N}$, with $\rho^0$ probability on $\R^{2d}$, $\rho^t$ is a solution to a {\it mean-field type (Vlasov)} {\bf nonlinear} equation of the following general form in the case of Section \ref{exach}

%\eqref{newever} in the case of Section \ref{exach} with  of the factorized form $\rho_N^0=(\rho^0)^{\otimes N}$, $\rho^t(x,v)$ satisfies the Vlasov equation

$$
\partial_t\rho^t(x,v)+v\partial_x\rho^t(x,v)=
$$
\be\label{metvla}
\partial_v\left(
\left(
\int
G(x,x';v,v')\rho(x',v')dx'dv'
+\eta
\nabla_{x}\varphi^t(x)
\right)
\rho(x,v)
\right)
\ee    
$$
\partial_s\psi^s(x)=\Delta\psi^s(x)
-\kappa\psi^s
%+\frac1N\sum_{j=1}^N\chi(x-x_i(t))
+\chi*\rho^s(x)\ \ \  0\leq s\leq t,
$$
with initial condition $\rho^0$.
\end{enumerate}
%\vfill
 The meaning of \eqref{meteuler} is easy to catch when no chemical interaction is present, i.e. when $\eta=0$. In this case the equation reduces to the equation satisfied by the push forward of the initial condition $\rho^0$ by the one particle dynamical system obtained by replacing the two body interaction $G$ by the average on the second of the two bodies by the solution itself (self-consistent field). It therefore reflect the fact that only one of the particle is kept moving, all the other being replaced by an average  pondered by the density itself.
 
 %{\hfill A-S. Sznitman $\sim$1990, \dots, F. Golse, C. Mouhot, P. 2017}
%\hfill (Cuccker-Smale)
\section{Propagation of chaos \ \small(recovering Independence)}
The tensorial power property of the initial condition $\rho_N^0=(\rho^0)^{\otimes N}$ is of fundamental importance for the derivation of the mean-field equation. Its meaning is that one starts with $N$ independent particles.

Of course this property immediately breaks down after the dynamics starts, i.e. $\rho_N^t\neq(\rho^t)^{\otimes N}$ for $t>0$ (except for trivial examples of dynamics): the dynamics couples particles and kills the independence property.

Nevertheless some trace of it remains valid when considering higher order marginals, namely $j$  particles probabilities of the form
\begin{eqnarray}
\rho^t
_{N;j}   
(x_1,\dots,x_j,v_1,\dots,v_j)&&\nonumber\\:=
\int\dot\int \rho^t_N(x_1,\dots,x_j,x_{j+1},\dots,x_N;v_1,\dots,v_j,v_{j+1},\dots,v_N)dx_2\dots dx_Ndv_2\dots dv_N&&\nonumber
\end{eqnarray}
for $j\leq N$ kept fixed.

Propagation of chaos states that, at any time $t>0$ and for any $j$,
$$
\lim_{N\to\infty}\rho^t
_{N;j}=(\rho^t)^{\otimes j}.
$$
\section{Velocity and density on (our) physical space \ \small(Euler)} 

Marginals again: in order to pass from a (single) one-particle function on phase-space, such as the solution $\rho^t$ of the Vlasov equation,  to two functions on physical space, we define
$$
\mu^t(x)=\int\rho^t(x,v)dv\hskip 3cm u^t(x)=\int v \rho^t(x,v)dv
$$
\hskip 5cm density\hskip 5cm velocity field    

\bigskip
Associated to these two quantities arises naturally  by induction from hydrodynamics the following Euler system of partial differential equations on physical space:
 
\be\label{meteuler}
\left\{
\begin{array}{l}
\partial_t\mu^t+\nabla(u^t \mu^t) = 0\\
\partial_t(\mu^tu^t)+
\nabla(\mu^t(u^t)^{ \otimes 2}) = 
%\mu^t
\int
f(\cdot-x';u^t(\cdot)-u^t(x')
\mu^t(x')dx'
%\right.\\
%\eta e^{-\tau t}\int_0^te^{(t-s)\Delta}\nabla\chi*\mu^s
%\left.
%\hskip 4.5cm+\ 
+\eta
\nabla\psi^t
\\
\partial_s\psi^s=\Delta\psi^s-\kappa\psi^s
+\chi*\rho^s
\end{array}
\right.
\ee 
A natural question then pops out of \eqref{meteuler}; wondering in which situations the density and velocity fields of a solution for the Vlasov equation \eqref{metvla}  are solution for the Euler system \eqref{meteuler}.

And indeed this is true:
\begin{itemize}
\item when $\rho^t(x,v)=\mu^t(x)\delta(v-u^t(x))$ monokinetic ansatz
\item for the case of the Krause model for agents systems (see below)
\item sometimes numerically, with both Cucker-Smale and chemotaxis (see Section \ref{ni} below).
\end{itemize}
\section{Density only on (our) physical space  \ \small(Keller-Segel)}\ 

The phenomenological Keller-Segel  describes the collective motion of cells by using  the density only, with no use of the velocity field (\cite{Pa, KS, Pe}). It reads

\be\label{KSS}
\left\{
\begin{array}{l}
\partial_t\mu^t-\Delta\mu^t-\nabla(\mu^t\varphi^t)=0\\
\alpha\partial_s\varphi^s=\Delta\varphi^s-\kappa\varphi^s+\mu^s,\ s\in[0,t].
%\\
%(u^0,\varphi^0)=(u^{in},\varphi^{in})\in H^s,\ s>\tfrac d2+1.
\end{array}
\right.
\ee
In \cite[Section 5.2]{Pe} was considered the case $\alpha=0$  and in \cite[Section 5.5.1]{Pe} the case where the Laplacian in \eqref{KSS} is replaced, e.g., by $\nabla(\mu^t\nabla\cdot)$. With these changes, the system \eqref{KSS} becomes
\be\label{KSSS}
\left\{
\begin{array}{l}
\partial_t\mu^t=\nabla(\mu^t\nabla\mu^t)+\nabla(\mu^t\varphi^t)\\
\Delta\varphi^t=\kappa\varphi^t-\mu^t.
%\\
%(u^0,\varphi^0)=(u^{in},\varphi^{in})\in H^s,\ s>\tfrac d2+1.
\end{array}
\right.
\ee
 
This system  can be derived from  a two-body chemotaxis system with friction of the form

\be\label{parteps}
\left\{
\begin{array}{rcl}
\dot xi(t)&=&\frac{v_i(t)}m\\
\dot v_i(t)&=&\frac{1}{N}\displaystyle\sum_{j=1}^N
\nabla V^\eps(x_i(t)-x_j(t))+\eta
\nabla_{x}\varphi^t(x_i(t))-v_i\\
\partial_t\varphi^t(x)&=&\Delta\varphi^t(x)
-\kappa\varphi^t
+\frac1N\sum\limits_{j=1}^N\chi^\eps(x-x_i(t)),\
\end{array}
\right.
\ee

Having in mind a large mass limit, we set $m=\eps^{-1}$.

A simple computation shows that, after a rescaling of time by $\eps$, the particle system \eqref{parteps} leads to
the following Euler system
\be\label{meteulereps}
\left\{
\begin{array}{l}
\partial_t\mu^t+\nabla(u^t \mu^t) = 0\\
\varepsilon\partial_t(\mu^tu^t)+
\nabla(\mu^t(u^t)^{ \otimes 2})= 
\mu^t
\int \nabla V^\eps(\cdot-y)\mu^t(y)dy
%f(\cdot-x';u^t(\cdot)-u^t(x')
%\mu^t(x')dx'
%\right.\\
%\eta e^{-\tau t}\int_0^te^{(t-s)\Delta}\nabla\chi*\mu^s
%\left.
%\hskip 4.5cm+\ 
+\eta\mu^t
\nabla\psi^t
-u^t\mu^t
\\
\eps\partial_s\psi^s=\Delta\psi^s-\kappa\psi^s
+\chi^\eps*\rho^s.
\end{array}
\right.
\ee 
%where the term $u^t\mu^t$ is a friction term (correspondig to add $v_i$ to the right-hand side of the last equation in \eqref{newever}, and the constant $\varepsilon$ can be seen as a joint rescaling of the mass $m_i$ in \eqref{newever}, time and the constant $\eta$.

Taking $V^\eps$ and $\chi^\eps$ as approximations of the Dirac mass and letting
%Indeed putting 
$\varepsilon\to 0$ in \eqref{meteulereps} gives
\be\label{umu}
u^t\mu^t=\eta\mu^t
\nabla\psi^t-+\mu^t\nabla \mu^t
\ee
and 
\be\label{umu2}
\Delta\psi^t=\kappa\psi^t
-\rho^t
\ee
so that we get  \eqref{KSSS} after  $u^t$ 
has been eliminated from the first equation in \eqref{meteulereps} thanks to \eqref{umu}, and having replaced the third equation in  \eqref{meteulereps} by \eqref{umu2}.
%, taking $V=\nabla\delta$ (or an approximation of it as $N\to\infty$) leads to \eqref{meteulereps}.

A rigorous proof of this derivation, obtained by  linking the limit $\eps\to0$ to $N\to\infty$ in \eqref{parteps}, together with other Keller-Segel type non-local equations corresponding to other asymptotic behaviors of the different parameters entering in \eqref{parteps}, can be found in \cite{TP}.

\part{THEOREMS\hfill {\large\it fini de rire !} {\scriptsize\it (J. Lacan)}}\label{theorems}

%\vskip 1cm
%\hfill {\large\it fini de rire !}
%
%\hfill {\scriptsize\it J. Lacan}
%\large

In this chapter, we will present rigorous results concerning the five examples already described.

As we  saw in Section \ref{examples}, all these examples are particular cases of the system
\eqref{newever}:
\be\nonumber
\left\{\begin{array}{rcl}
\dot x_i(t)&=&\frac{v_i(t)}{m_i},\  m_i\in (0,+\infty],\  i=1,\dots,N\\
\dot v_i(t)&=&F_i\left(t,x_1(t),\dots,x_N(t),v_1(t),\dots,v_N(t),\{x_j(s),j=1,\dots,N,0\leq s<t\}\right).
\end{array}
\right.
\ee  

In the sequel we will denote by $\rho_{N;s}, s=0,\dots,N$ the $s$-marginal of a probability $\rho$ on $\R^{2dN}$ defined by
\be\label{defmar}
\rho_{N,s}(x_1,\dots,x_s,v_1,\dots,v_s)
:=
\int_{\R^{2d(N-s}}\rho(x_1,\dots,x_N,v_1,\dots,v_N)dx_{s+1}dv_{s+1}\dots dx_Ndv_N)
\ee
with the convention $\rho_{N;0}=\rho$ and $\rho_{N,N}=1$.

Let us recall that we denote in Section \ref{micprob} the mapping $\Phi_N^t$ which associates to any initial condition the corresponding solution of \eqref{newever}.
\vskip 0.5cm
In the sequel we will use the two Wasserstein distances defined as follow (\cite{V1,V2}).

\begin{Def}[Wasserstein distances]\label{defwas}\ 

For all probability measures $\mu_1,\mu_2$ on $\bR^m$
$$
W_1(\mu_1,\mu_2)=\sup\left\{\int_{\R^m}f\,d(\mu_1-\mu_2)\ \mid\ f\in\Lip(\R^m),\ \Lip(f)\leq 1\right\}.
$$

For all probability measures $\mu,\nu$ on $\bR^m$ with finite second moments 
$$
\MKd(\mu,\nu)^2
=
\inf_{\gamma\in\Gamma(\mu,\nu)}\int_{\bR^m\times\bR^m}
|x-y|^2\gamma(dx,dy)
$$
where $\Gamma(\mu,\nu)$ is the set of probability measures on $\bR^m\times\bR^m$ whose marginals on the two factors are $\mu$ and $\nu$.
\end{Def}
\vskip 1cm
%Let us recall that all examples trated in this paper have the following form:
%
%\be\nonumber
%\left\{\begin{array}{rcl}
%\dot x_i(t)&=&\frac{v_i(t)}{m_i},\  m_i\in (0,+\infty],\  i=1,\dots,N\\
%\dot v_i(t)&=&
%\frac1N\sum\limits_{j=1}^NF^t\left(\left\{\left(x_i(t)-x_j(s),v_i(t)-v_j(s)),
%%1\leq j\leq N,
%\ 0\leq s\leq t
%%\right),\right.\right.\\ 
%%&&\left.\left.j=1,\dots,N,\ 0\leq s\leq t\}
%\right)\right\}\right)\\
%F^t(\dots)&=&
%F(x_i(t)-x_j(t))+\eta G^t_{\{
%x_i(t)-x_j(s),\ 0\leq s\leq t\}}(x_i(t))+\mu F^{\mbox{\scriptsize ext}}(x_i(t))\\
%\partial_t G^t_{\{
%x-x_j(s), 0\leq s\leq t\}}(x)&=&HG^t_{\{
%x-x_j(s),\ 0\leq s\leq t\}}(x)+J(x-x_j(s)).
%\end{array}
%\right.
%\ee  
%with $\eta=0$ in the nonchemiotaxis examples, and 
%$$
%H=\Delta+\tau,\ J(y)=\nabla\chi(y)
%$$
%in the chemiotaxis case.
%
%In the general framework, the Vlasov equation reads
In the sequel
\begin{itemize}
\item $\Phi^t_N$ will be the mapping defined through \eqref{newever}:
 
\begin{eqnarray}
\Phi^t_N(x_1(0),\dots,x_N(0)t;v_1(0),\dots,v_N(0))&&\nonumber\\
:=(x_1(t),\dots,x_N(t);v_1(t),\dots,v_N(t))
 \mbox{ solution to \eqref{newever}}&&\nonumber
 \end{eqnarray}
 \item $\rho^t$ will be solution of the Vlasov equation \eqref{metvla} with initial condition 
 $\rho^{in}$
 \item $\psi^t$ ill be the potential part of the solution to the Vlasov equation \eqref{metvla} with initial condition $\psi^{in}=0$.
\end{itemize}

\section{Bounded two-body (Hamiltonian) interactions}\label{btbi}\ 

In this case, all $m_i=1$ and $F_i(\dots)=-\tfrac1N\sum\limits_{j=1}^N\nabla V(x_i-x_j),\ \nabla V\mbox{ Lipschiz}$.

The Vlasov equation in this case is the ``standard"  one and reads

\be\label{vlsbtb}
\partial_t\rho(x,v)=-v\cdot\nabla_x\rho(x,v)
+\int\nabla V(x-y)\rho(y,w)dydw\cdot \nabla_v\rho(x,v).
\ee
The result is
\begin{Thm}[\cite{GMP}]\ 

For each  $j$,
$$
W_2\left((\Phi^t_N\#(\rho^{in})^{\otimes N})_{N;j},(\rho^t)^{\otimes j}\right)
\leq e^{\Lambda t}N^{-\frac12}
$$
for some $\alpha,\Lambda>0$.
\end{Thm}

\section{Cucker-Smale type systems}\label {cssts}
%\end{document}

%\newcommand{\bR}{\mathbb R}
%\newcommand{\gen}{\gamma}
%\newcommand{\vla}{\nu}
%\newcommand{\vlau}{\mathbf v}
%\newcommand{\ta}{\tau}
%\newcommand{\hist}[2]{[#1]^{\leq #2}}
%\newcommand{\histk}[2]{[#1]_k^{\leq #2}}
%
%\newcommand{\vlav}[3]{\vlau(#1,#2,\hist{#3}{#1})}
%\newcommand{\vlavk}[3]{\vlau(#1,#2,\histk{#3}{#1})}
%
%
%
%\newcommand{\vlaun}{\mathbf \Gamma}
%\newcommand{\histn}[2]{[(#1)^{\otimes N}]^{\leq #2}}
%\newcommand{\histkn}[2]{[(#1)^{\times N}]_k^{\leq #2}}
%\newcommand{\vlavn}[3]{\vlaun(#1,#2,(\histn{#3}{#1})_{N;1})}
%\newcommand{\vlavkn}[3]{\vlaun^k(#1,#2,(\histkn{#3}{#1)_{N:1}})}
%\newcommand{\vlavkni}[3]{\vlaun^k_i(#1,#2,(\histkn{#3}{#1)_{N:1}})}

%\newcommand{\rhoi}{\rho^{in}}
%\newcommand{\rhoin}{(\rhoi)^{\otimes N}}
%
%\newcommand{\macphi}{\varphi}
%\newcommand{\macphin}{\varphi^{in}}
%
%\newcommand{\rc}{{\tau_c}}

%\newcommand{\cblue}[1]{\begin{color}{blue}#1\end{color}}
%
%\newcommand{\rod}{\tau}

\def\cA {\mathcal{A}}
\def\cB {\mathcal{B}}
\def\cC {\mathcal{C}}
\def\cD {\mathcal{D}}
\def\cE {\mathcal{E}}
\def\cF {\mathcal{F}}
\def\cG {\mathcal{G}}
\def\\fH {\mathcal{H}}
\def\cJ {\mathcal{J}}
\def\cK {\mathcal{K}}
\def\cL {\mathcal{L}}
\def\cM {\mathcal{M}}
\def\cN {\mathcal{N}}
\def\cP {\mathcal{P}}
\def\cQ {\mathcal{Q}}
\def\cR {\mathcal{R}}
\def\cS {\mathcal{S}}
\def\cT {\mathcal{T}}
\def\cU {\mathcal{U}}
\def\cV {\mathcal{V}}
\def\cW {\mathcal{W}}
\def\cX {\mathcal{X}}
\def\cZ {\mathcal{Z}}

%\newcommand{\vvec}{\overrightarrow}
%\newcommand{\Div}{\operatorname{div}}
%\newcommand{\Rot}{\operatorname{curl}}
%\newcommand{\Diam}{\operatorname{diam}}
%\newcommand{\Dom}{\operatorname{Dom}}
%\newcommand{\Sign}{\operatorname{sign}}
%\newcommand{\Span}{\operatorname{span}}
%\newcommand{\Supp}{\operatorname{supp}}
%\newcommand{\Det}{\operatorname{det}}
%%\newcommand{\Tr}{\operatorname{trace}}
%\newcommand{\Codim}{\operatorname{codim}}
%\newcommand{\Dist}{\operatorname{dist}}
%\newcommand{\DDist}{\operatorname{Dist}}
%\newcommand{\Ker}{\operatorname{Ker}}
%\newcommand{\IM}{\operatorname{Im}}
%\newcommand{\Id}{\operatorname{Id}}
%\newcommand{\Lip}{\operatorname{Lip}}
%\newcommand{\Osc}{\operatorname{osc}}
%\newcommand{\Esssup}{\operatorname{ess sup}}
%\newcommand{\Essinf}{\operatorname{ess inf}}
%\newcommand{\ad}{\operatorname{\mathbf{ad}}}
%\newcommand{\AD}{\operatorname{\mathbf{AD}}}
%\newcommand{\conj}{\operatorname{\mathbf{conj}}}
%\newcommand{\wtost}{\operatorname{wto}^*}
%\newcommand{\Subsetc}{\,\mathop{\subset}_c\,}

% \newcommand{\MKd}{W_2}
%%%%%%%%%%%%%%%%%%%%%%%%%%%%%%%%%%%%%%%%%%%%%%%%%%%%%%%%%%%%%

In this case all the $m_i$ are equal to $1$ and
\be\label{defG}
%F_i(t,Y,W)=\frac1N\sum_{j=1}^N\gamma(w_i-w_j,y_i-y_j)+\eta
%\nabla_{z}%\varphi(t,z,Y)|_{z=y_i}
F_i(t,X,V)=\frac1N\sum_{j=1}^N\gamma(v_i-v_j,x_i-x_j)
\ee
%$X(t)=(x_1(t),\dots,x_N(t))$, 
$\gamma$ is the collective interaction function supposed to be Lipschitz continuous.

The Vlasov equation, first stated in \cite{HT} reads

\be\label{vlscs}
\partial_t\rho(x,v)=-v\cdot\nabla_x\rho(x,v)
+\nabla_v\cdot\left(\int\gamma(x-y,v-w)\rho(y,w)dydw\rho(x,v)\right).
\ee
Note that when $\gamma(x,v)=\nabla V(x)$ \eqref{vlscs} reduces to \eqref{vlsbtb}.

The result holds true.
\begin{Thm}[\cite{NP}]\ 

For each  $j$,
$$
W_2\left((\Phi^t_N\#(\rho^{in})^{\otimes N})_{N;j},(\rho^t)^{\otimes j}\right)
\leq C_tN^{-\alpha}
$$
for some $\alpha>0$.
\end{Thm}

\section{Topological models}\label{topmod}
\newcommand{\fcr}{\mathcal X}
In this case one takes $m_i=1,i=1,\dots,N$ and
\begin{equation}
  \label{eq:pij}
  F_i(\dots)=\tfrac1N\sum_{j=1}^N K\bigl(M(x_i,|x_i - x_j|)\bigr)(v_i-v_j),
\end{equation}
where $K
\colon [0,1] \to \R^+$ is a positive Lipschitz continuous non-increasing
function,
%\sout{such that $\int_0^1 K(z) \, \de z= \gamma$}
and, for $r>0$, the rank function
\begin{equation}
  \label{eq:Mpart}
  M(x_i,r) = \frac 1N \sum_{k=1}^N \fcr\{|x_k - x_i| \le  r\}
\end{equation}
counts the number of agents at a distance
less than or equal to $r$ from $x_i$, normalized with $N$.
Note that in this case $p_{ij}$ is a stepwise function of the positions
of all the agents and $F_i$ is not of type ``two body" anymore, it depends on all agents.

The Vlasov equation was introduced in \cite{BCR} and  reads
\be\label{vltop}
\partial_t\rho(x,v)=-v\cdot\nabla_x\rho(x,v)
+
\ee

Note that it is not quadratic in $\rho$ anymore.

The result is the following.
\begin{Thm}[\cite{BPR}]\ 
For each  $j$,
$$
W_1\left((\Phi^t_N\#(\rho^{in})^{\otimes N})_{N;j},(\rho^t)^{\otimes j}\right)
\leq C_tN^{-\alpha}
$$
for some $\alpha>0$.
\end{Thm}

\section{chemotaxis}\label{sctx}
In this case all the $m_i$ are equal to $1$ and
\be\label{defG}
%F_i(t,Y,W)=\frac1N\sum_{j=1}^N\gamma(w_i-w_j,y_i-y_j)+\eta
%\nabla_{z}%\varphi(t,z,Y)|_{z=y_i}
F_i(t,X,V)=\frac1N\sum_{j=1}^N\gamma(v_i-v_j,x_i-x_j)+\eta
\nabla_{x}\varphi^t(x_i)
%+F_{ext}(y_i),
+F_{ext}(x_i),
\ee
%$X(t)=(x_1(t),\dots,x_N(t))$, 
$\gamma$ is the collective interaction function supposed to be Lipschitz continuous, $F_{ext}$ is an external force  and $\varphi$ satisfies the equation
\be\label{defeqphi0}
%\partial_s\varphi(s,z, Y(t))=D\Delta_z\varphi-\kappa\varphi +f(z,Y(s)),
\partial_s\varphi^s(x)=D\Delta_x\varphi-\kappa\varphi +f(x,X(s)),
\ s\in[0,t],\ \varphi^{s=0}=\macphin
%\in C^2(\bR^d),
\ee 
for some $\kappa, D,\eta\geq0$ and function $f$ of the form
\be\label{deff}
f(x,X)=\frac1N\sum_{j=1}^N\chi(x-x_i),\ \ \chi\in\mathcal C_c^1.
\ee
%In \eqref{eq1}, \eqref{defeqphi0} and \eqref{deff}, $X=(x_1,\dots,x_N)\ V=(v_1,\dots,v_N)$.

%Moreover t
%The  function $\gamma:\bR^d\times \bR^d\to\bR^d$ is supposed to be Lipschitz continuous.

\vskip 1cm

%{\bf Particle system}: $Z^{in}:(X_N^{in},V_N^{in})\to %(X_N(t),V_N(t))=\Phi_N^t(X_N^{in},V_N^{in})
%Z^t=(X^t,V^t)=\Phi_N^t(Z^{in})
%$
%\begin{eqnarray}
%\dot x_i=v_i,\  \dot v_i(t)&=&\frac{1}{N}\displaystyle\sum_{j=1}^N
%F(x_i(t)-x_j(t),v_i(t)-v_j(t))+\eta
%\nabla_{x}\varphi^t(x_i(t))
%\nonumber\\
%\partial_t\varphi^t(x)&=&\Delta\varphi(x)
%%-\kappa\varphi 
%+\frac1N\sum_{j=1}^N\chi(x-x_i(t)),
%\nonumber
%\end{eqnarray}
%\centerline{$\chi:\mbox{ Lipschitz supported ion a ball {\LARGE\color{red}$\bullet$}}.
%$}
 
The Vlasov equation, in fact a system when considering the potential $\varphi$ introduced in \cite{NP0} reads
%\scriptsize
%$$
%\partial_t\rho^t(x,v)+v\partial_x\rho^t(x,v)
%$$
%$$=
%\partial_v\left(
%\left(
%\int
%F(x-x';v-v')\rho(x',v')dx'dv'
%+\eta
%\nabla_{x}\varphi^t(x)
%\right)
%\rho(x,v)
%\right)
%$$
%$$
%\partial_s\psi^t(x)=\Delta\psi(x)
%%-\kappa\varphi 
%+\chi *\rho^s,\ \ \ 0\leq s\leq t,\ \ \ \psi^{s=0}=0
%$$

\be\label{vlche}
\left\{
\begin{array}{ll}
\partial_t\rho^t(x,v)+v\partial_x\rho^t(x,v)
&\\
=\partial_v\left(
\left(
\int
F(x-x';v-v')\rho(x',v')dx'dv'
+\eta
\nabla_{x}\varphi^t(x)
\right)
\rho(x,v)
\right)
&\\
\rho^{t=0}:=\rho^{in} \\
\partial_s\psi^t(x)=\Delta\psi(x)
%-\kappa\varphi 
+\chi *\rho^s,\ \ \ 0\leq s\leq t,\ \ \ \psi^{s=0}=0
\end{array}
\right.
\ee
\vskip 1cm
In order to state the result in this case, we need to fix the notation. 
 
 With the notation $Z^{in}:(X_N^{in},V_N^{in})\to %(X_N(t),V_N(t))=\Phi_N^t(X_N^{in},V_N^{in})
Z^t=(X^t,V^t)=\Phi_N^t(Z^{in}),\ X^t=(x_1(t),\cdot,x_N(t))
$ and the same for  $V$ and $Z$,

The {\it non-Markovian} particle system reads then

\begin{eqnarray}\label{vphi}
\dot x_i=v_i,\  \dot v_i(t)&=&\frac{1}{N}\displaystyle\sum_{j=1}^N
F(x_i(t)-x_j(t),v_i(t)-v_j(t))+\eta
\nabla_{x}\varphi^t(x_i(t))
\nonumber\\
Z^{t=0}&:=&Z^{in}\nonumber\\
\partial_t\varphi^t_{Z^{on}}(x)&=&\Delta\varphi^t_{Z^{on}}(x)
%-\kappa\varphi 
+\frac1N\sum_{j=1}^N\chi(x-x_i(t)),
\end{eqnarray}

%\centerline{$\chi:\mbox{ Lipschitz supported ion a ball {\LARGE\color{red}$\bullet$}}.
%$}

since  the solution of \eqref{vphi} at time $t$ depends on the whole trajectory $\{X^s,\ 0\leq s\leq t\}$ that is on $Z^{in}$. By the same argument, the last equation of the Vlasov system 
\eqref{vlche} should be stated this way:
$$
\partial_s\psi^t_{\rho^{in}}(x)=\Delta\psi^t_{\rho^{in}}(x)
%-\kappa\varphi 
+\chi *\rho^s,\ \ \ 0\leq s\leq t,\ \ \ \psi^{s=0}=0
$$

The question which naturally arises is then: how to compare $\varphi^t_{Z^{in}}$ (particles) which depends on $Z^{in}$
to $\psi^t_{\rhoi}$ (Vlasov) which depends on $\rhoi$? 

The answer is given in our next result.

\begin{Thm}[\cite{NP0}]\ 

 Let $\varphi^t_{Z^{in}}$ the chemical potential of the particle system starting at $X^{in}$,
  and let $\psi^t_{\rhoi}$ the chemical potential solution to the Vlasov system \eqref{vlche}.  
  
{Then
$$
W_2\left((\Phi^t_N\#(\rho^{in})^{\otimes N})_{N;1},\rho^t\right)^2\leq 
%C(t)
%N^{-\alpha}.
%\sqrt{C_d(N)}
\tau(t)
%(\alpha(t)^2+{\beta(t)^2+\gamma(t)^2})
%\times
%%\left\{
%%\begin{array}{ll}
%%N^{-\frac12}&d=1\\
%%N^{-\frac12}{\log{N}}&d=2\\
%%N^{-\frac1{d}}&d>2
%\end{array}
%\right.
N^{-\tfrac1{d}}\ \ \  \mbox{\scriptsize $(d>2)$}
$$
and    

%{\it \color{black} How to compare $\varphi^t_{Z^{in}}$ (particles) which depends on $Z^{in}$
%
%\hfill to $\psi^t_{\rhoi}$ (Vlasov) which depends on $\rhoi$?}    
$$
\int_{\bR^{2dN}}\|\nabla\varphi^t_{Z^{in}}-\nabla\psi^t_{\rhoi}\|^2_\infty
(\rhoi)^{\otimes N}(dZ^{in})
%\leq t^2\Lip(\nabla\chi)^2\MKd\left((\Phi^t_N\#(\rho^{in})^{\otimes N})_{N;1},\rho^t\right)^2
\leq 
%C(t)
%N^{-\alpha}.
%\sqrt{C_d(N)}
%t^2\Lip(\nabla\chi)^2
\rc(t)
%5t^2\Lip(\nabla\chi)^2(\tau(t)+Ce^{\Gamma(t)}
%%C
%)
%(\alpha(t)^2+{\beta(t)^2+\gamma(t)^2})
%\left\{
%\begin{array}{ll}
%N^{-\frac12}&d=1\\
%N^{-\frac12}{\log{N}}&d=2\\
N^{-\tfrac1{d}}\ \ \  \mbox{\scriptsize $(d>2)$}
%\end{array}
%\right.
$$}
 \end{Thm}
 
 \section{Hydrodynamics limit}\label{hydlim}\LARGE
 
In the case of Vlasov evolution with a monokinteic initial condition $\mu^{in}(x)\delta(v-u^{in}(x))$, the following result holds true in the chemotaxis case.
\begin{Thm}[\cite{NP0}]\  

%\small
Let {$\mu^t,u^t$ solves the Euler system:}
$$
\left\{
\begin{array}{l}
\partial_t\mu^t+\nabla(u^t \mu^t) = 0\\
\partial_t(\mu^tu^t)+
\nabla(\mu^t(u^t)^{ \otimes 2}) = 
%\mu^t
\int
f(\cdot-x';u^t(\cdot)-u^t(x')
\mu^t(x')dx'
%\right.\\
%\eta e^{-\tau t}\int_0^te^{(t-s)\Delta}\nabla\chi*\mu^s
%\left.
%\hskip 4.5cm+\ 
+\eta
\nabla\psi^t
\\
\partial_s\psi^s=\Delta\psi^s-\kappa\psi^s
+\psi^s*\rho^s
\end{array}
\right.
$$
%\vskip 1cm
Then {$\mu^t(x)\delta(v-u^t(x))$ solves the Vlasov system \eqref{vlche}} with initial condition $\mu^{in}(x)\delta(v-u^{in}(x))$.
\end{Thm}
A similar result is valid in the topological interactions type situation, see \cite{BPR}.
%\section{Hydrodynamiocs limit}\label{hydlim}
%\Large\newpage\Large
\LARGE
\section{Numerical intermezzo \cite{MNP}}\label{ni}\ 
In this review we have considered three different ways of handling large number of particle systems: the ODE original (and true) vision \eqref{newever}, its mean-field limit Vlasov equation \eqref{metvla} and its hydrodynmaics picture Euler equation \eqref{meteuler}.

%\vskip 3cm
%\hfill\it\Large to press or not to press
Let us briefly show in this section some numerical results (therefore done for a finite number $N$ of particles) showing first that the kinetic (Vlasov) description  is quite faithful to the original particle one, and second how the hydrodynamical picture (Euler) can be compared to the kinetic one.
\subsection{Particles $\longrightarrow$ Vlasov}\ 
%\scriptsize

Here we compare the particle evolution  with the corresponding Vlasov one for three different times. the two first simulations correspond to a flocking situation (with zero momentum of the center of gravity), and the third one a non-flocking one.

\LARGE
\begin{figure}[h!]
	\centering
	%se vuoi le lettere sotto le sottofigure, mettere [][], e nel secondo slot mettere cosa è quella sottofigura.Tipo:
	%\subfigure[][Initial configuration]{\includegraphics[width=6cm, height=2.3cm]{figs/Vlas_b005_t0}}
	 particles 
	 
	 \subfigure{\includegraphics[width=3cm, height=3cm]{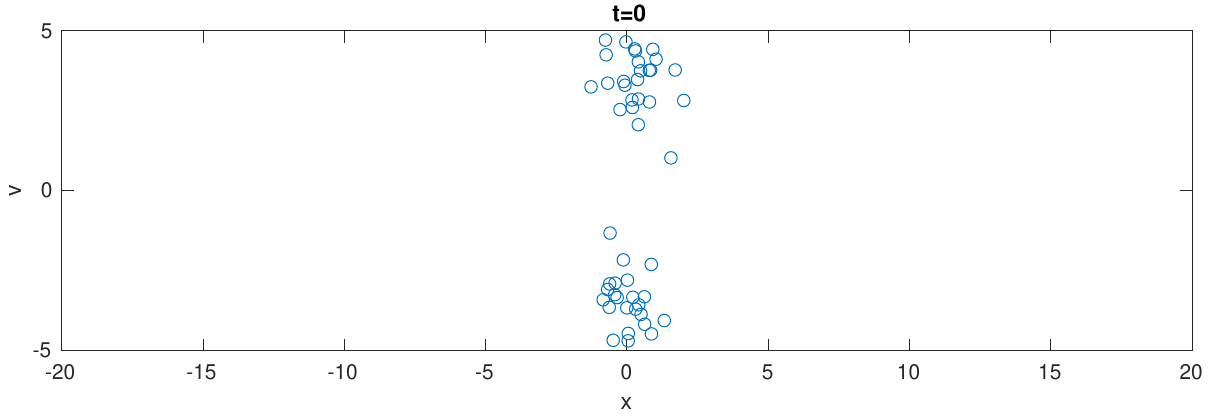}}
	\hspace{0 mm}
	\subfigure{\includegraphics[width=3cm, height=3cm]{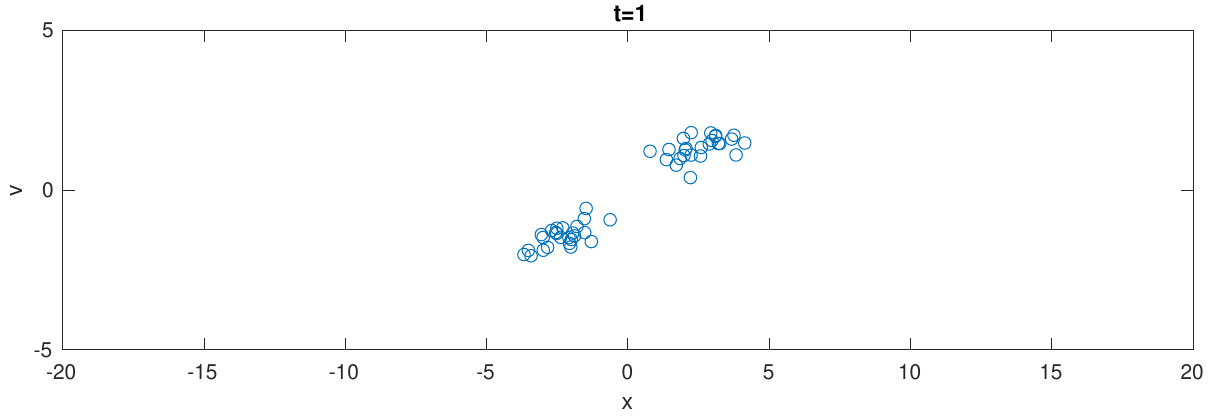}}
	\hspace{0 mm}
	\subfigure{\includegraphics[width=3cm, height=3cm]{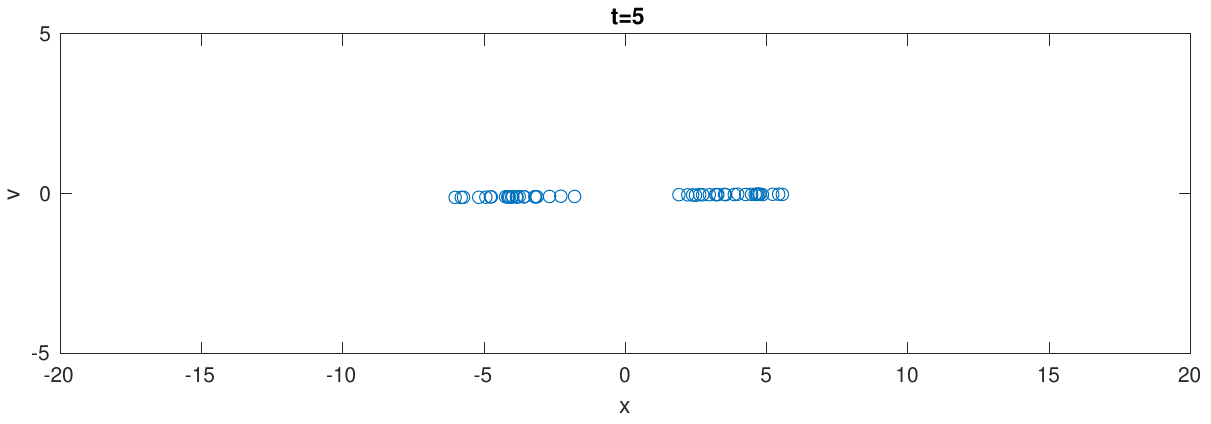}}\\
	kinetic 
	
	\subfigure{\includegraphics[width=3cm, height=3cm]{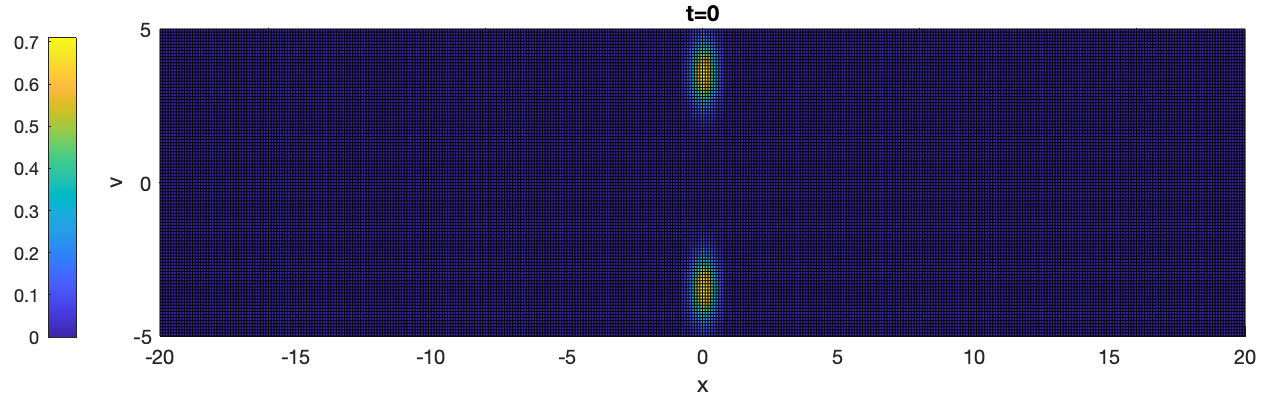}}
	\hspace{0 mm}
	\subfigure{\includegraphics[width=3cm, height=3cm]{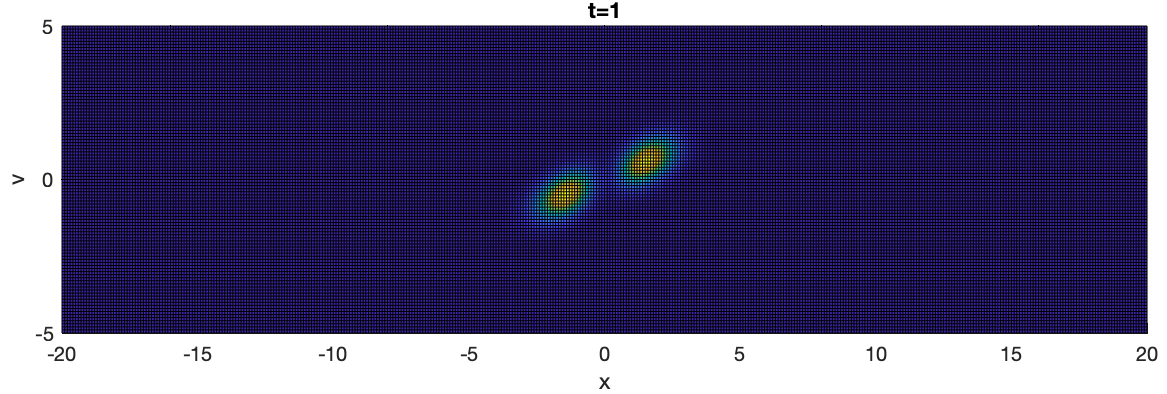}}
	\hspace{0 mm}
	\subfigure{\includegraphics[width=3cm, height=3cm]{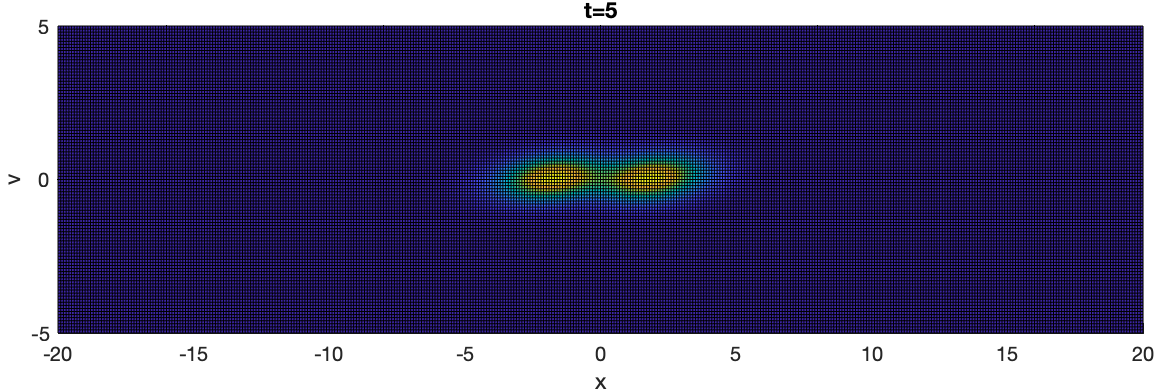}}
	%\caption{Test 1: Numerical simulation of Cucker-Smale model with $\beta=0.05$, at particle level (first line) and kinetic (second line) level.}
	\label{fig:Vlasov_b005}
\end{figure}

\begin{figure}[h!]
	\centering
	%se vuoi le lettere sotto le sottofigure, mettere [][], e nel secondo slot mettere cosa è quella sottofigura.Tipo:
	%\subfigure[][Initial configuration]{\includegraphics[width=6cm, height=2.3cm]{figs/Vlas_b005_t0}}
%	\hspace{-13 mm}
%	\subfigure{\includegraphics[width=10.9cm, height=3.4cm]{figs/Vlas_b005e095_t0}}
%	\hspace{0 mm}
%	\subfigure{\includegraphics[width=10cm, height=3.3cm]{figs/Vlaschemo_b095_t0dot6}}
%	\hspace{0 mm}
%	\subfigure{\includegraphics[width=10cm, height=3.3cm]{figs/Vlaschemo_b095_t5}}
\hspace{-13 mm}
particles

	\subfigure{\includegraphics[width=3cm, height=3cm]{figs/CS_b095_t0}}
	\hspace{0 mm}
	\subfigure{\includegraphics[width=3cm, height=3cm]{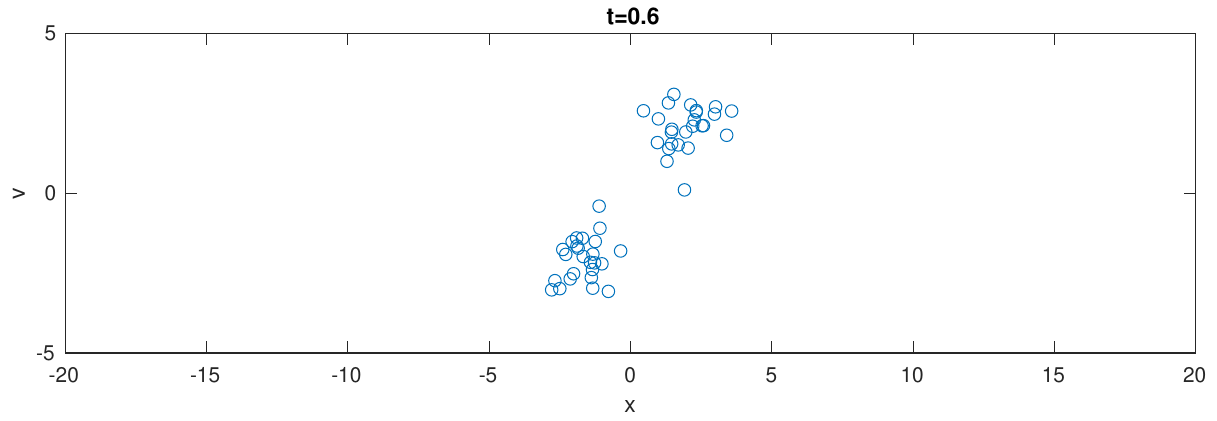}}
	\hspace{0 mm}
	\subfigure{\includegraphics[width=3cm, height=3cm]{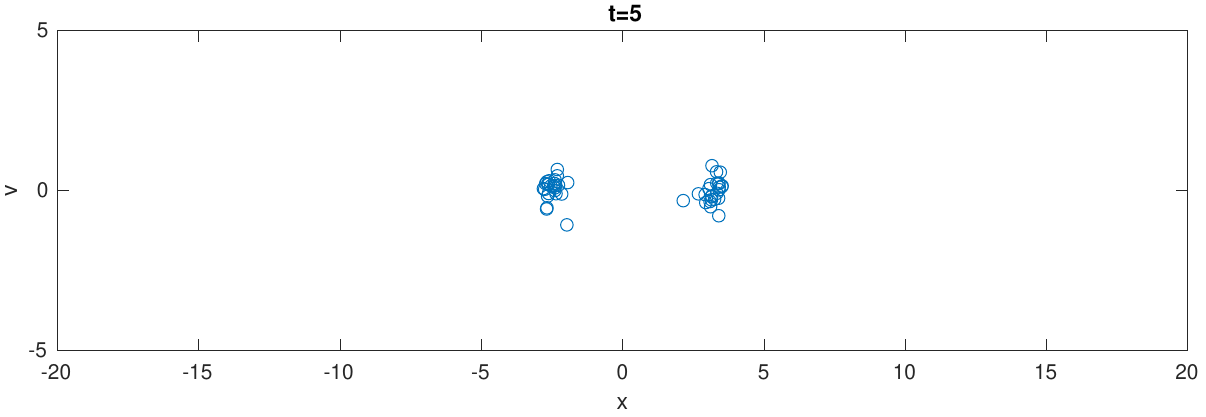}}\\
kinetic 

	\hspace{-13 mm}
	\subfigure{\includegraphics[width=3cm, height=3cm]{figs/Vlas_b005e095_t0}}
	\hspace{0 mm}
	\subfigure{\includegraphics[width=3cm, height=3cm]{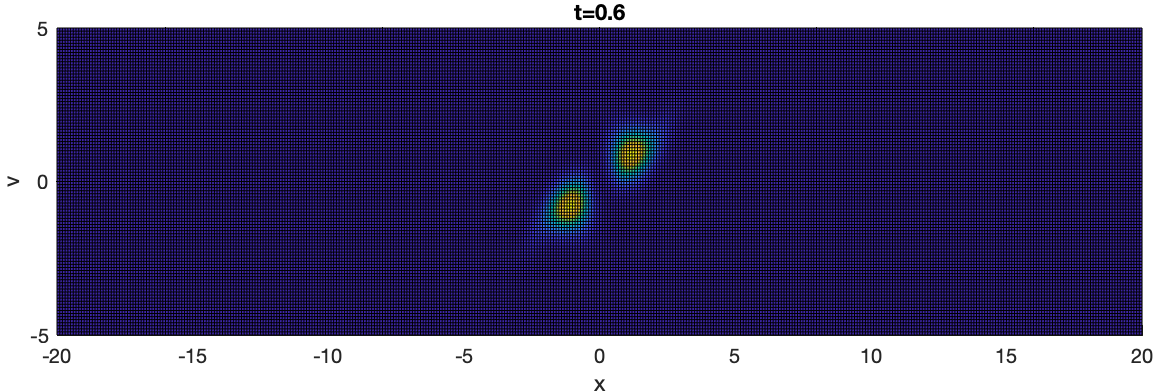}}
	\hspace{0 mm}
	\subfigure{\includegraphics[width=3cm, height=3cm]{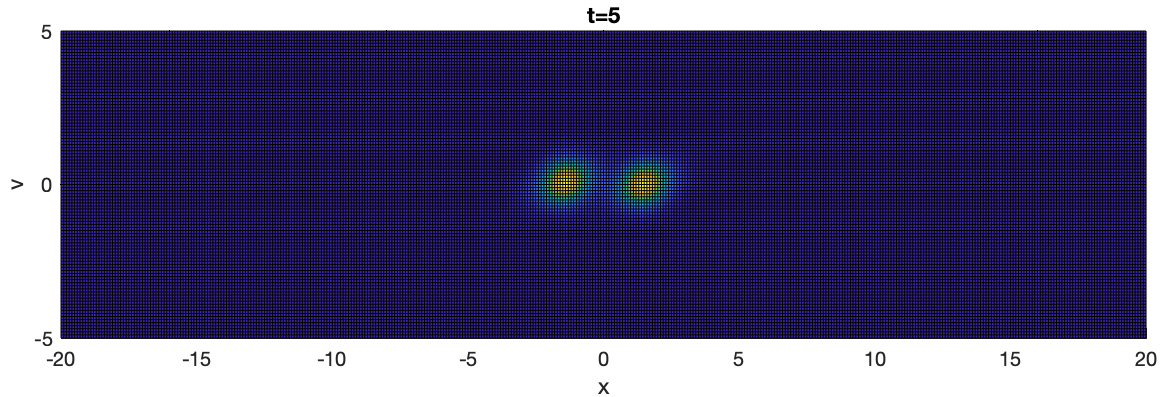}}
	%\caption{Test 3: Numerical simulation of Cucker-Smale model with chemotaxis at particle level (first line) and kinetic (second line) level, with parameters as in Test 2 and $\eta=1.4$.}
	\label{fig:Vlasov_b095_chemo}
	\end{figure}

	\begin{figure}[h!]
	\centering
	%se vuoi le lettere sotto le sottofigure, mettere [][], e nel secondo slot mettere cosa è quella sottofigura.Tipo:
	%\subfigure[][Initial configuration]{\includegraphics[width=6cm, height=2.3cm]{figs/Vlas_b005_t0}}
%	\hspace{-13 mm}
%	\subfigure{\includegraphics[width=10.9cm, height=3.4cm]{figs/t0}}
%	\hspace{0 mm}
%	\subfigure{\includegraphics[width=10cm, height=3.3cm]{figs/t2}}
%	\hspace{0 mm}
%	\subfigure{\includegraphics[width=10cm, height=3.3cm]{figs/t3}}
	\hspace{-13 mm}
	particles 
	
	\subfigure{\includegraphics[width=3cm, height=3cm]{figs/CS_b095_t0}}
	\hspace{0 mm}
	\subfigure{\includegraphics[width=3cm, height=3cm]{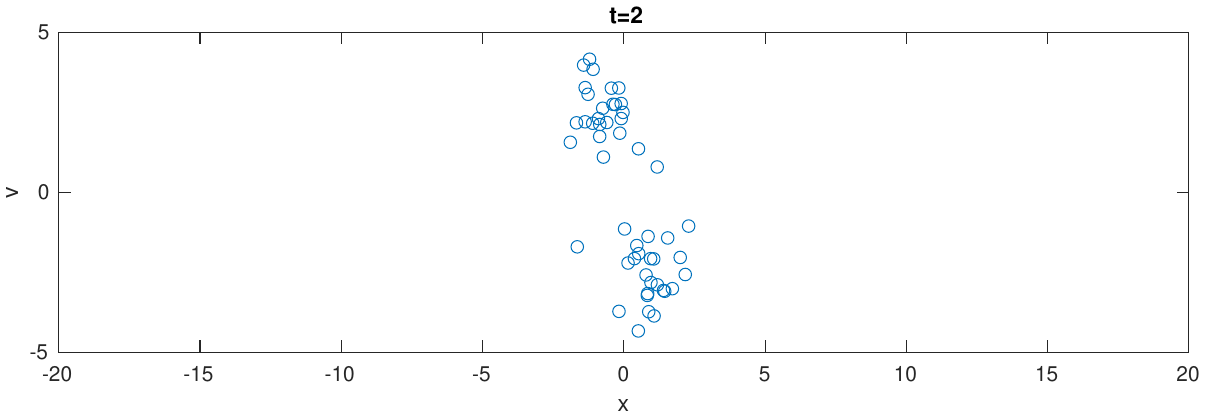}}
	\hspace{0 mm}
	\subfigure{\includegraphics[width=3cm, height=3cm]{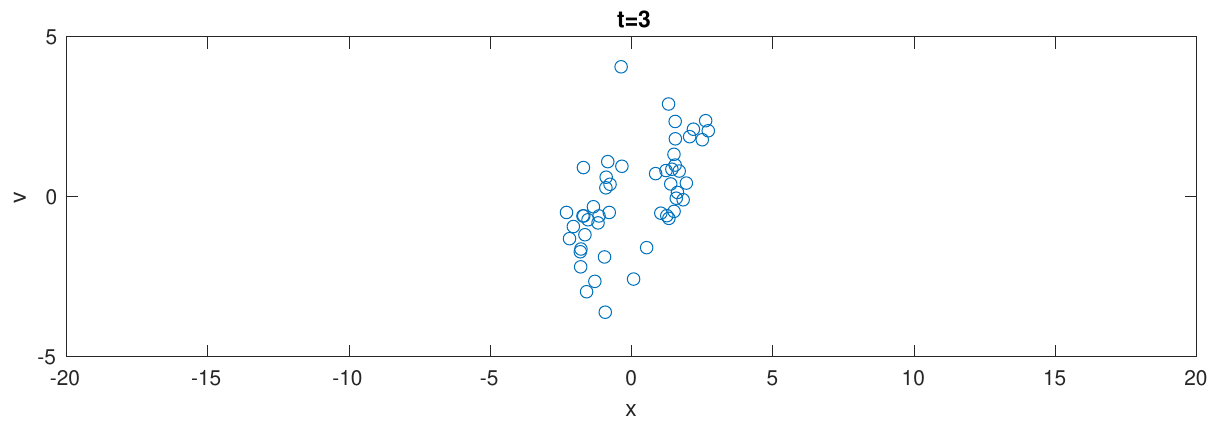}}\\
	\hspace{-13 mm}
	
	kinetic 
	
	\subfigure{\includegraphics[width=3cm, height=3cm]{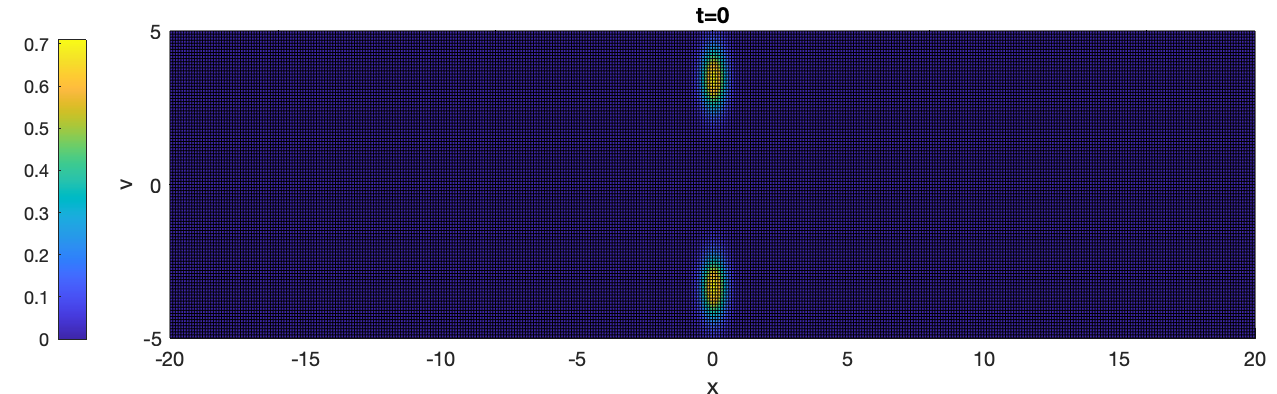}}
	\hspace{0 mm}
	\subfigure{\includegraphics[width=3cm, height=3cm]{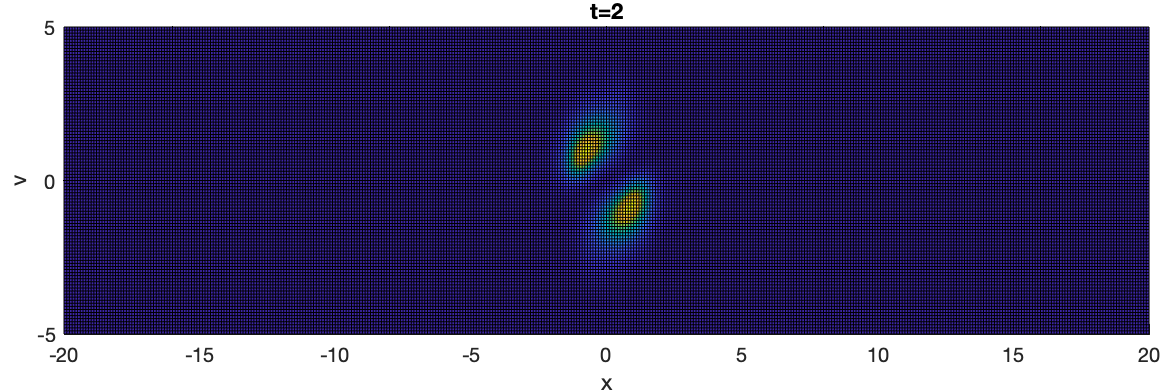}}
	\hspace{0 mm}
	\subfigure{\includegraphics[width=3cm, height=3cm]{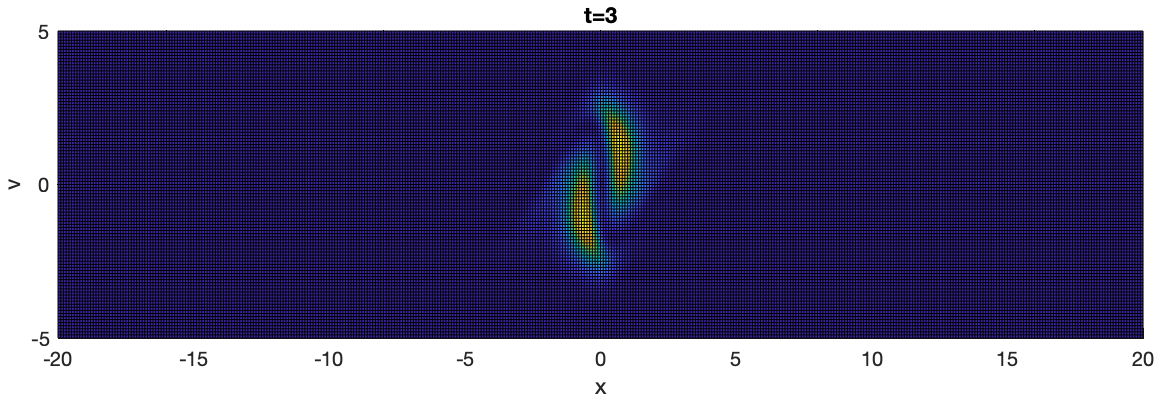}}
	%\caption{Test 4:
	\label{fig:Vlasov_onlychemo}
\end{figure}
\newpage
%\centerline{\small particles $\longrightarrow$ Vlasov}

 \subsection{Vlasov $\longrightarrow$ Euler}\ 
 
 Let us compare now the evolution through Vlasov to the one of the density and velocity fields through Euler. More precisely, we compare the density and the velocity field of the solution to the Vlasov equation (blue curve) with the evolution through Euler equation of the density and velocity field of the same initial condition to Vlasov (dashed curves).
 
 The first line shows the initial condition, the second and third ones the evolution (at time $t=2$)  without chemotaxis and  with chemotaxis respectively.

\begin{figure}[h!]
	\centering
%\raisebox{10pt}{\parbox[b]{1.2\textwidth}{ \hspace{2cm} Vlasov phase space \hspace{1.5cm} $\nu_0^t$ and $\mu^t$ \hspace{2.5cm} $\nu_1^t$ and $Q^t$}}\\
%\raisebox{10pt}{\parbox[b]{1.2\textwidth}{ \hspace{3cm} \small{Vlasov phase space} \hspace{2.5cm} $\nu_0^t$ and $\mu^t$ \hspace{3.5cm} $\nu_1^t$ and $Q^t$}}\\
%\raisebox{70pt}{\parbox[b]{.12\textwidth}{$t=2$\\\\$\eta=0$\\ $\alpha=0$}}

\hskip 1.1cm Vlasov\hskip 1.3 cm denstity/Euler\hskip 1cm velocity/Euler

\subfigure{\includegraphics[width=3cm, height=3cm]{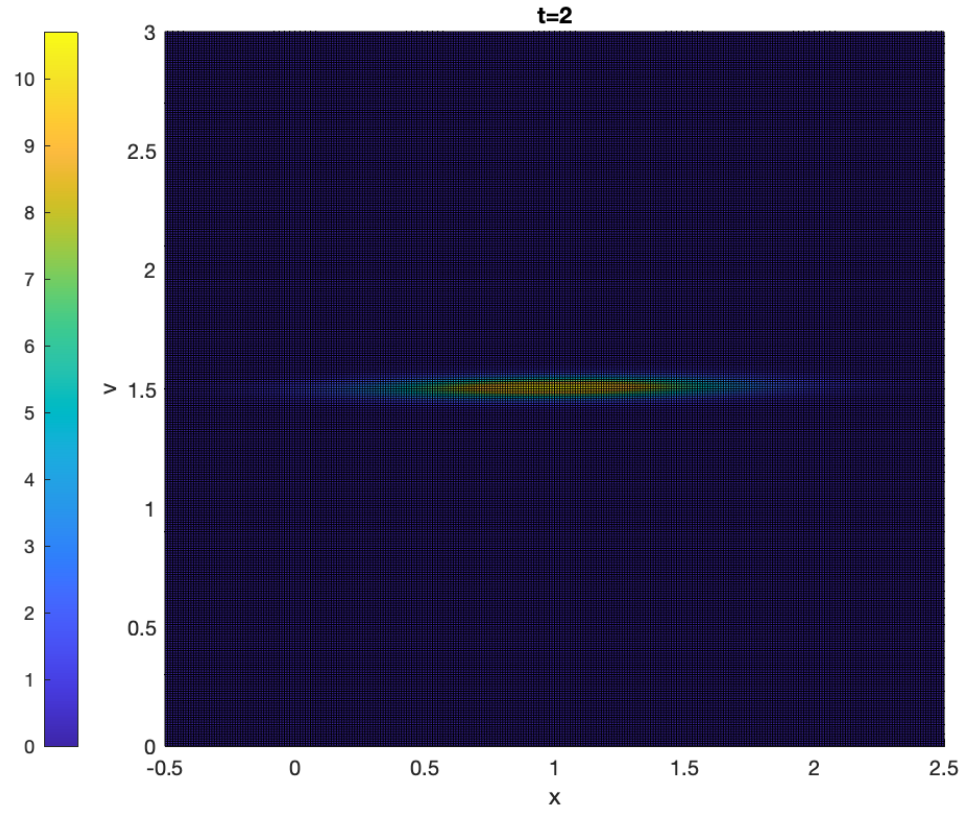}}\hskip 1cm
\subfigure{\includegraphics[width=3cm, height=3cm]{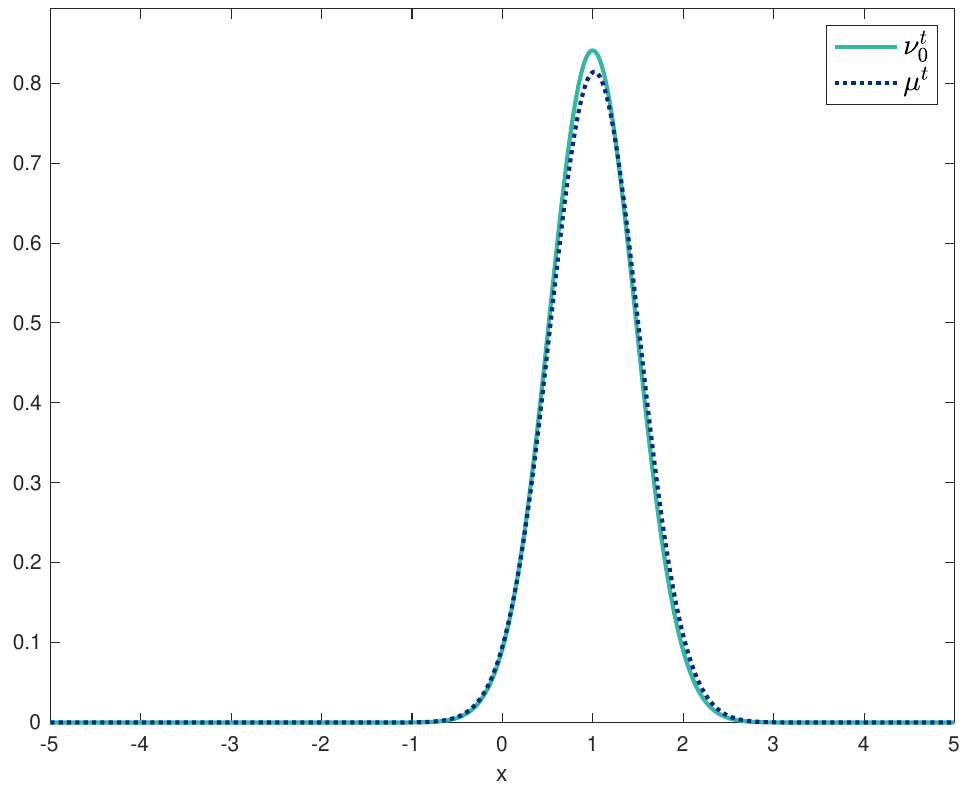}}\hskip 1cm
\subfigure{\includegraphics[width=3cm, height=3cm]{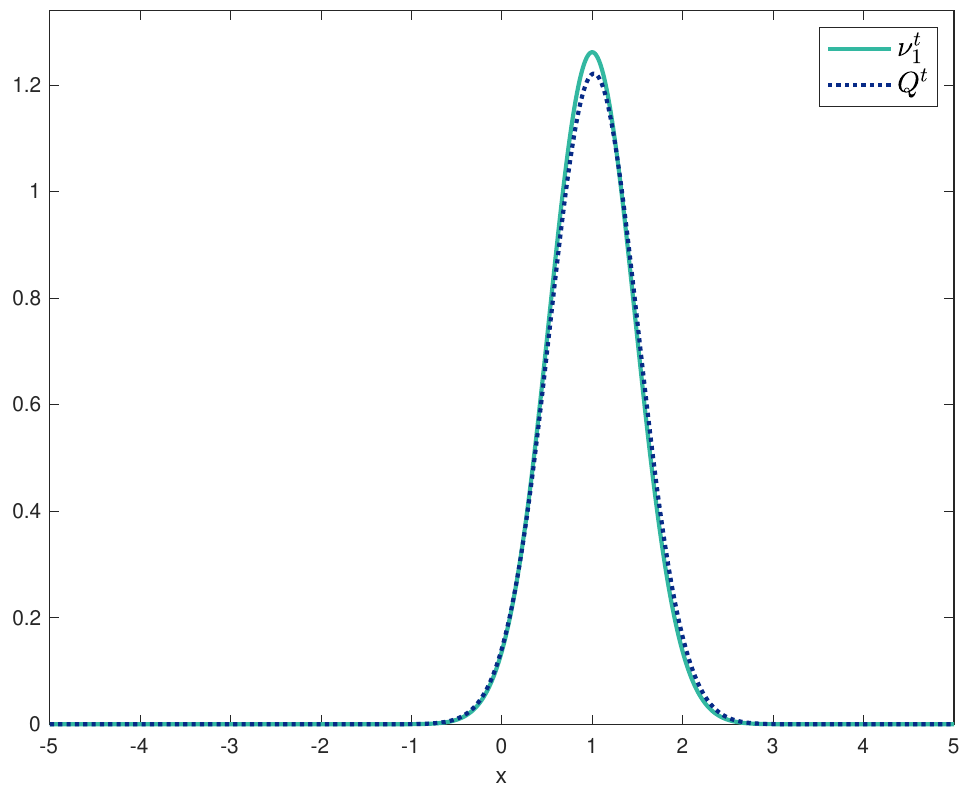}}\\
%\raisebox{70pt}{\parbox[b]{.12\textwidth}{$t=2$\\\\$\eta=0.2$ \\ $\alpha=0$}}
\hspace{0.5cm}
\subfigure{\includegraphics[width=3cm, height=3cm]{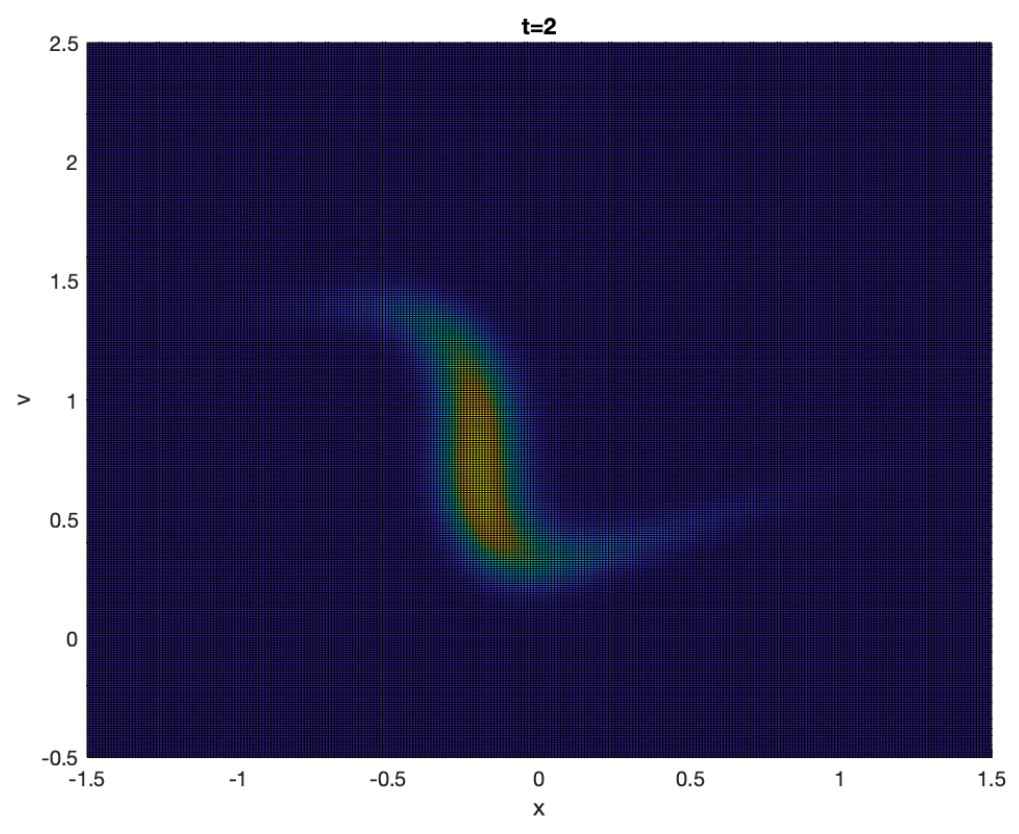}}\hskip 1cm
\subfigure{\includegraphics[width=3cm, height=3cm]{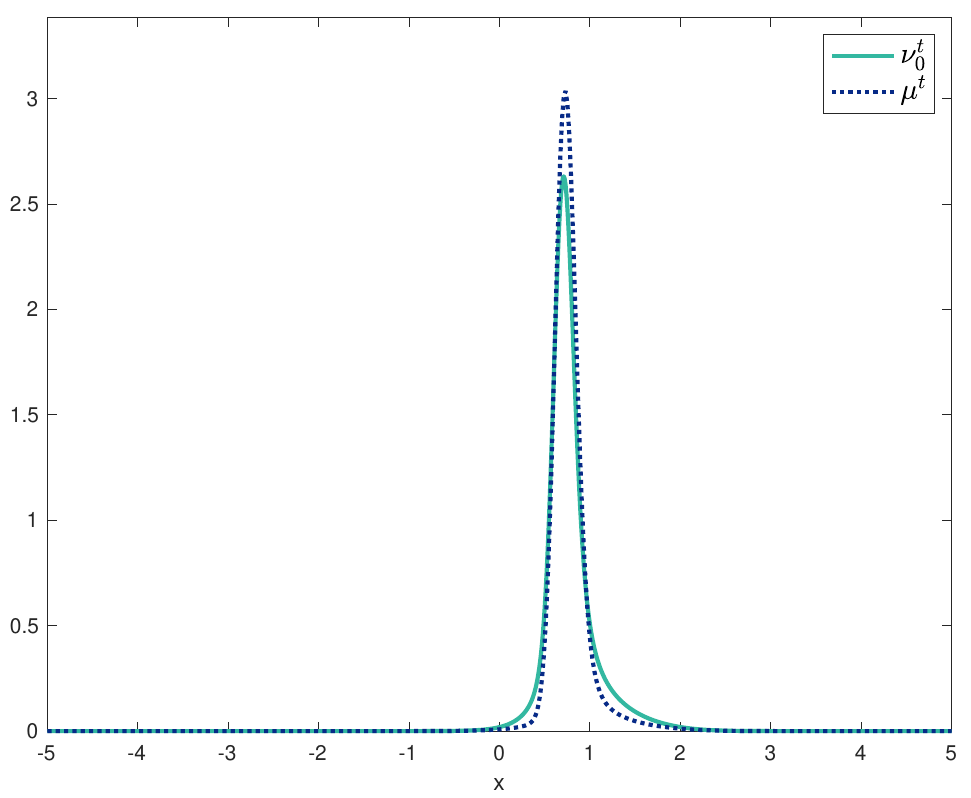}}\hskip 1cm
\subfigure{\includegraphics[width=3cm, height=3cm]{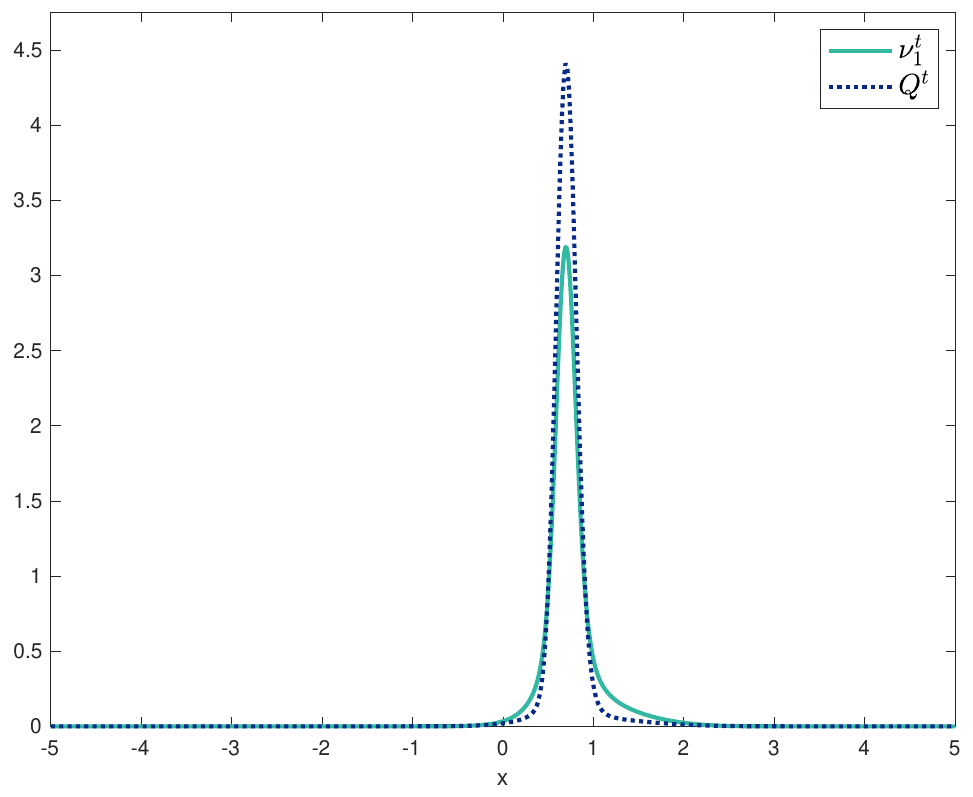}}	\\
%\raisebox{70pt}{\parbox[b]{.12\textwidth}{$t=3$\\\\$\eta=0.2$ \\ $\alpha=2$}}
\hspace{0.5cm}
\subfigure{\includegraphics[width=3cm, height=3cm]{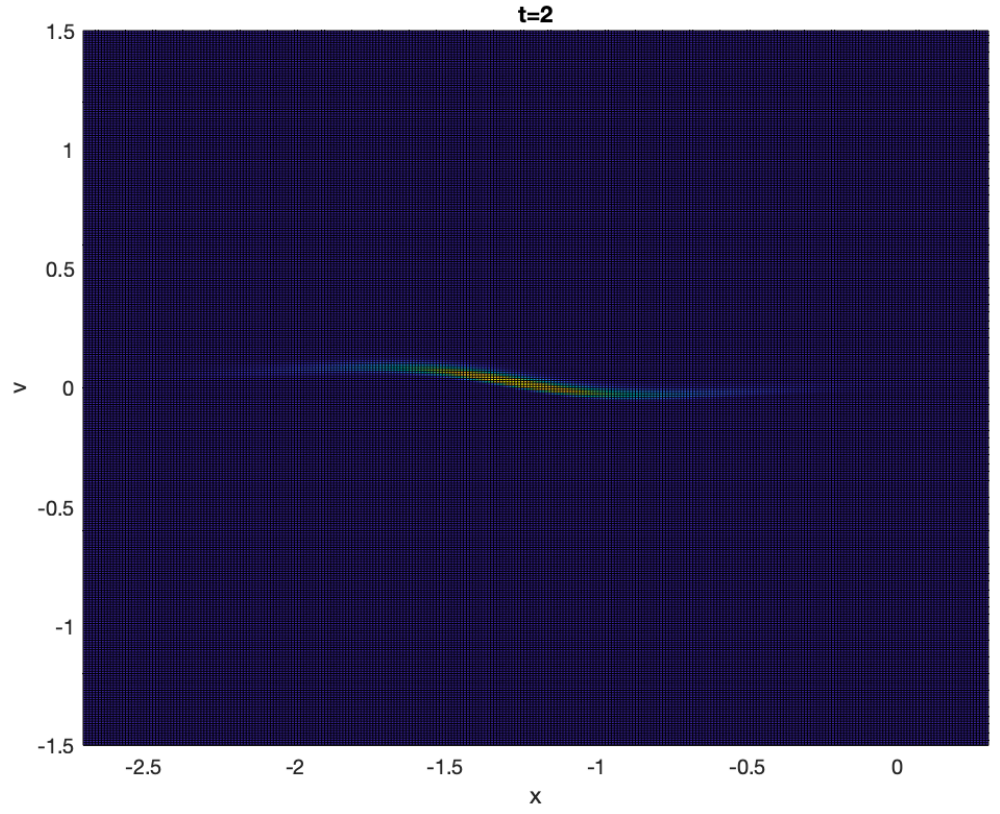}}\hskip 1cm
\subfigure{\includegraphics[width=3cm, height=3cm]{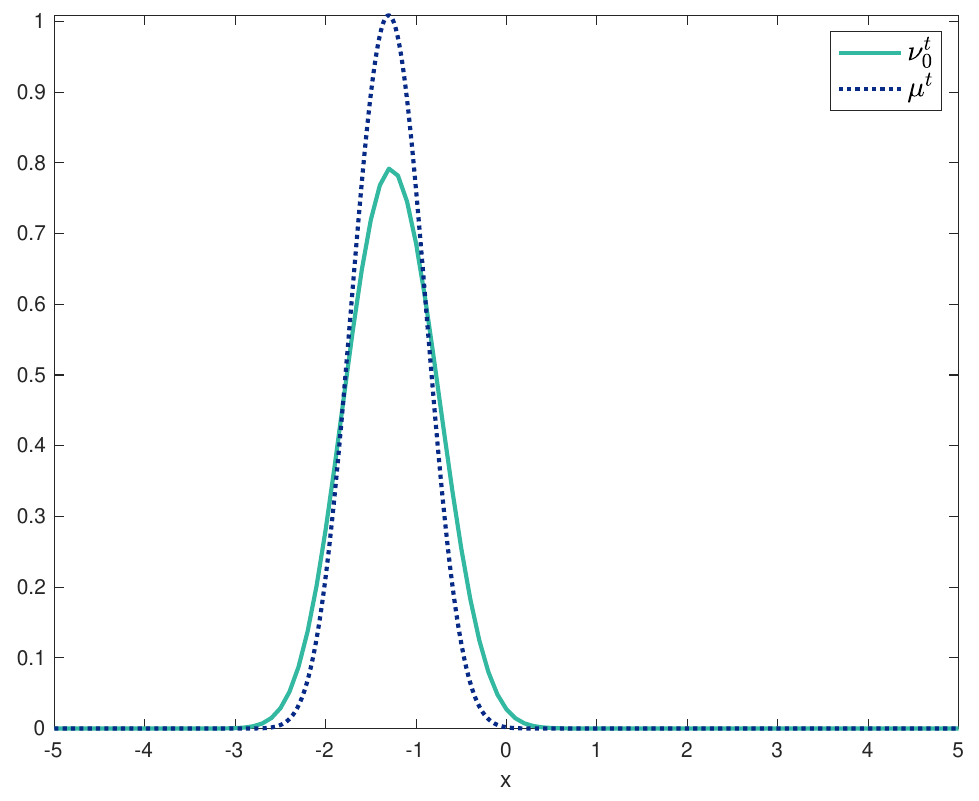}}\hskip 1cm
\subfigure{\includegraphics[width=3cm, height=3cm]{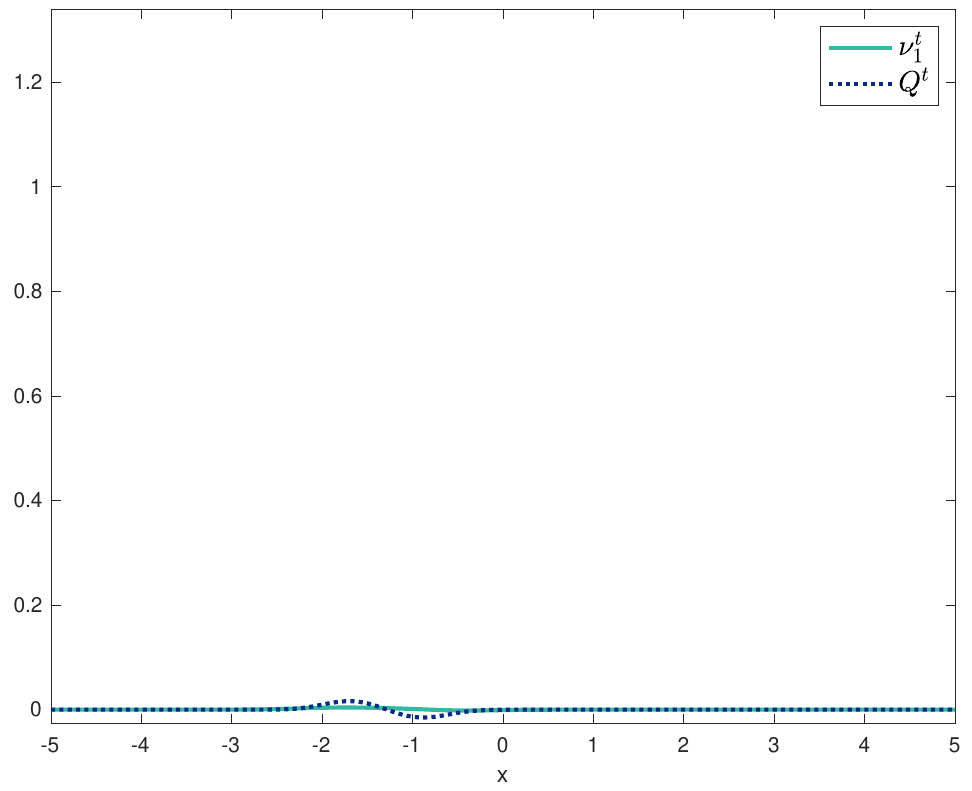}}	
%\subfigure{\includegraphics[width=5cm, height=5cm]{figs/density_c05_eta1}}
%\subfigure{\includegraphics[width=5cm, height=5cm]{figs/velocity_c05_eta1}}
	%\caption{Test 7: approximation of monokinetic initial data. Comparison between zero-th order moment of the solution of (V) and $\mu^t$ (second column),  first order moment of the solution of (V) and $Q^t$ (third column): without chemotaxis (first line: $\eta=0$, $\alpha=0$ at $t=2$), with chemotaxis (second line: $\eta=0.2$, $\alpha=0$ at $t=2$), with chemotaxis and damping (third line: $\eta=0.2$, $\alpha=2$ at $t=2$).}
	\label{fig:monokinetic}
\end{figure}
% {\scriptsize 
%\hfill M. Menci, R. Natalini, T. P.
%
%\hfill Microscopic, kinetic and hydrodynamic hybrid models of collective motions with
%
%\hfill chemotaxis : a numerical study (2023)}
\subsection{Moral}
Instead of commenting each figure, let us exhibit the main features of these numerical computations (see  \cite{MNP} for (much) more details):
\begin{itemize}
\item in many
situations the most striking features of the finite particle level  are preserved by passing
 to the one of Vlasov
 \item the fidelity of Euler versus Vlasov
is increased by the presence of the chemical interaction
\item and this even without monokineticity hypothesis (therefore outside of the paradigm of strict equivalence Vlasov $\sim$ Euler
\item the nonlocal integral Euler system  keeps memory of the interactions at
the microscopic level 
\item adding an additional pressure term of size $\epsilon$, still keeping the nonlocal
integral term but with no equivalent term in Vlasov, helps
\item there exists of an  {optimal value $\epsilon$}  realizing an improved correspondence between Vlasov moments and Euler solutions
\end{itemize}

\LARGE

\end{document}